\documentclass[titlepage,draft,12pt]{article} 
\usepackage{amsfonts,latexsym,pstricks,pst-plot} 
\usepackage[a4paper]{geometry}

\newcommand{\ocont}{$\omega$-continuous}

\newcommand{\aholder}{$\alpha$-H\"{o}lder}
\newcommand{\G}{\mathcal{G}}
\newcommand{\trebar}[1]{\|#1\|}
\newcommand{\dainfty}{D(A^{\infty})}

\newcommand{\lk}{\lambda_{k}}
\newcommand{\dk}{\delta_{k}}
\newcommand{\sqrdk}{\sqrt{\delta_{k}}}

\newcommand{\ep}{\varepsilon}
\newcommand{\re}{{\mathbb{R}}}
\newcommand{\n}{{\mathbb{N}}}

\newcommand{\vs}{\vspace{.2cm}}

\newcommand{\prf}{{\sc Proof.}$\;$}
\newcommand{\qed}{{\penalty 10000\mbox{$\quad\Box$}}\bigskip}

\newcommand{\foralll}{\forall\:}

\newcommand{\m}[1]{m(\au{#1}^{2})}
\newcommand{\au}[1]{|A^{1/2}#1|}
\newcommand{\auq}[1]{|A^{1/2}#1|^{2}}

\newcommand{\cep}{c_{\ep}}

\newcommand{\uep}{u_{\ep}}

\newtheorem{thm}{Theorem}[section]
\newtheorem{thmbibl}{Theorem}
\newtheorem{rmk}[thm]{Remark}
\newtheorem{prop}[thm]{Proposition}

\newtheorem{cor}[thm]{Corollary}
\newtheorem{ex}[thm]{Example}
\newtheorem{lemma}[thm]{Lemma}

\title{Derivative loss for Kirchhoff equations with non-Lipschitz
nonlinear term}
\author{Marina Ghisi\vs\\ {\normalsize
Universit\`a degli Studi di Pisa} \\{\normalsize Dipartimento di
Matematica ``Leonida Tonelli''}\\
{\normalsize 
PISA (Italy)}\\  
{\normalsize e-mail: \texttt{ghisi@dm.unipi.it}}\and 
Massimo Gobbino\vs\\ {\normalsize Universit\`a degli Studi di Pisa} 
\\{\normalsize Dipartimento di Matematica Applicata ``Ulisse Dini''}\\ 
{\normalsize 
 PISA (Italy)}\\  
{\normalsize e-mail: \texttt{m.gobbino@dma.unipi.it}}}
\date{}

\begin{document}
\maketitle
\begin{abstract}
	In this paper we consider the Cauchy boundary value problem for
	the integro-differential equation
	$$u_{tt}-m\left(\int_{\Omega}^{}|\nabla u|^{2}\,dx\right)
	\Delta u=0 
	\hspace{3em}
	\mbox{in }\Omega\times[0,T)$$
	with a continuous nonlinearity $m:[0,+\infty)\to[0,+\infty)$.
	
	It is well known that a \emph{local} solution exists provided that
	the initial data are regular enough.  The required regularity
	depends on the continuity modulus of $m$.
	
	In this paper we present some counterexamples in order to show
	that the regularity required in the existence results is sharp, at
	least if we want solutions with the same space regularity of
	initial data.  In these examples we construct indeed local
	solutions which are regular at $t=0$, but exhibit an instantaneous
	(often infinite) derivative loss in the space variables.
	
\vspace{1cm}

\noindent{\bf Mathematics Subject Classification 2000 (MSC2000):}
35L70, 35L80, 35L90.

\vspace{1cm} 
\end{abstract}
 
\section{Introduction}

In this paper we consider the hyperbolic partial differential equation
\begin{equation}
	u_{tt}(t,x)-
	m{\left(\int_{\Omega}\left|\nabla u(t,x)\right|^2\,dx\right)}
	\Delta u(t,x)=0
	\hspace{2em}
	\forall(x,t)\in\Omega\times[0,T),
	\label{eq:k}
\end{equation}
where $\Omega\subseteq\re^{n}$ is an open set, $\nabla u$ and $\Delta 
u$ denote the gradient and the Laplacian of $u$ with respect to the
space variables, and $m:[0,+\infty)\to[0,+\infty)$. Equation
(\ref{eq:k}) is usually considered with initial conditions 
\begin{equation}
	u(x,0)=u_{0}(x),
	\hspace{2em}
	u_{t}(x,0)=u_{1}(x)
	\hspace{2em}
	\forall x\in\Omega,
	\label{eq:data}
\end{equation}
and boundary conditions, for example of Dirichlet type (but the theory
is more or less the same also with Neumann or periodic boundary
conditions)
\begin{equation}
	u(x,t)=0
	\hspace{2em}
	\forall(x,t)\in\partial\Omega\times[0,T).
	\label{eq:boundary}
\end{equation}

From the mathematical point of view, (\ref{eq:k}) is probably the
simplest example of quasilinear hyperbolic equation.  From the mechanical
point of view, this Cauchy boundary value problem is a model for the
small transversal vibrations of an elastic string ($n=1$) or membrane
($n=2$).  In the string context it was introduced by 
\textsc{G.\ Kirchhoff} in \cite{kirchhoff}.

A lot of papers have been devoted to existence of local or global
solutions to (\ref{eq:k}), (\ref{eq:data}), (\ref{eq:boundary}). The
results can be divided into four main families.
\begin{enumerate}
	\renewcommand{\labelenumi}{(\Alph{enumi})} 
	\item \emph{Local existence in Sobolev spaces}.  If $m$ is a
	locally Lipschitz continuous function such that
	$m(\sigma)\geq\nu>0$ for every $\sigma\geq 0$ (strict
	hyperbolicity), then (\ref{eq:k}), (\ref{eq:data}),
	(\ref{eq:boundary}) admits a unique local solution for initial
	data in Sobolev spaces.  This result was first proved by
	\textsc{S.\ Bernstein} in the pioneering paper~\cite{bernstein},
	and then extended with increasing generality by many authors.  The
	more general statement is probably contained in \textsc{A.\
	Arosio} and \textsc{S.\ Panizzi}~\cite{ap} (see also the
	references quoted therein), where it is proved that the problem is
	locally well posed in the phase space
	$$V_{\beta}:=H^{\beta+1}(\Omega)\times H^{\beta}(\Omega)\ +\
	\mbox{boundary conditions}$$
	for every $\beta\geq 1/2$.

	\item \emph{Global existence for analytic data}.  If $m$ is a
	continuous function such that $m(\sigma)\geq 0$ for every
	$\sigma\geq 0$ (weak hyperbolicity), and $u_{0}(x)$ and $u_{1}(x)$
	are analytic functions, then (\ref{eq:k}), (\ref{eq:data}),
	(\ref{eq:boundary}) admits at least one global solution
	($T=+\infty$), which is analytic in the space variables for any
	time.
	
	This result was proved with increasing generality by \textsc{S.\
	Bernstein}~\cite{bernstein}, \textsc{S. I.
	Pohozaev}~\cite{poho-an}, \textsc{A.\ Arosio} and \textsc{S.\
	Spagnolo}~\cite{as}, \textsc{P.\ D'Ancona} and \textsc{S.\
	Spa\-gno\-lo}~\cite{das-an-1,das-an-2}.
	
	\item \emph{Local existence with assumptions between (A) and (B)}.
	More recently, \textsc{H.\ Hirosawa} \cite{hirosawa-main}
	considered the continuity modulus $\omega$ of the nonlinear term
	$m$.  He proved existence of local solutions in suitable classes
	of initial data, depending on $\omega$, both in the strictly
	hyperbolic and in the weakly hyperbolic case.  From the point of
	view of local solutions these results represent an interpolation
	between the results of (A) and (B).  The precise relations between
	$\omega$ and the regularity required on initial data are stated in
	section~\ref{sec:statements}.  Roughly speaking, in the strictly
	hyperbolic case the situation is summed up by the following
	scheme:
		\begin{eqnarray*}
			\omega(\sigma)=o(1) & \to & \mbox{analytic data,}  \\
			\omega(\sigma)=\sigma^{\alpha}\mbox{ (with $\alpha\in(0,1)$)} & \to
			& \mbox{Gevrey space
			$\G_{s}(\Omega)$ with $s=(1-\alpha)^{-1}$,}  \\
			\omega(\sigma)=\sigma|\log \sigma| & \to & 
			\mbox{$C^{\infty}(\Omega)$ (or $V_{\beta}$
			with finite derivative loss),}  \\
			\omega(\sigma)=\sigma & \to & \mbox{$V_{1/2}$ (with no
			derivative loss)}.
		\end{eqnarray*}
	
	More regularity is required in the weakly hyperbolic case, as
	shown in the following scheme: 
		\begin{eqnarray*}
			\omega(\sigma)=o(1) & \to & \mbox{analytic data,}  \\
			\omega(\sigma)=\sigma^{\alpha}\mbox{ (with $\alpha\in(0,1)$)} & \to
			& \mbox{Gevrey space
			$\G_{s}(\Omega)$ with $s=1+\alpha/2$,}  \\
			\omega(\sigma)=\sigma & \to & \mbox{Gevrey space
			$\G_{3/2}(\Omega)$}.
		\end{eqnarray*}
	
	\item \emph{Global existence in special situations}.  Besides (B),
	there are four special cases in which global solutions are known
	to exist.  We refer the interested reader to the quoted literature
	for the details.  We just point out that in all these results the
	nonlinearity $m$ is assumed to be Lipschitz continuous and
	strictly positive.
	\begin{itemize}
		\item  \textsc{K.\ Nishihara}~\cite{nishihara} proved global
		existence for quasi-analytic initial data. This class of
		functions strictly contains the space of analytic functions,
		but it is strictly contained in every Gevrey class
		$\G_{s}(\Omega)$ with $s>1$.
	
		\item \textsc{J.\ M.\ Greenberg} and \textsc{S.\ C.\
		Hu}~\cite{gh}, and then \textsc{P.\ D'ancona} and \textsc{S.\
		Spagnolo}~\cite{das} proved global existence for small initial
		data in Sobolev spaces in the case $\Omega=\re^{n}$, where
		dispersion plays a crucial role.  Later on this
		dispersion-based approach was extended to external domains
		(see \cite{yamazaki1,yamazaki2} and the references quoted
		therein).
	
		\item  \textsc{S.\ I.\ Pohozaev}~\cite{poho-m} proved global
		existence for initial data in the Sobolev space $V_{1}$
		in the case where $m(\sigma)=(a+b\sigma)^{-2}$, with $a>0$,
		$b>0$. In this particular case indeed equation (\ref{eq:k})
		admits a second order invariant.
	
		\item In some recent papers \textsc{R.\
		Manfrin}~\cite{manfrin1} (see also \cite{manfrin2},
		\cite{hirosawa2}, \cite{gg:k-manfrinosawa}) proved global
		existence in a new class of nonregular initial data.
		Manfrin's spaces are small in the sense that they don't
		contain any Gevrey class $\mathcal{G}_{s}(\Omega)$ with $s>1$,
		but they are large in the sense that any initial condition
		$(u_{0},u_{1})\in V_{1}$ is the sum of two initial conditions
		belonging to these spaces!
	\end{itemize}
\end{enumerate}

Despite of the many positive results, as far as we know no paper has
been devoted to negative results.  In particular in the literature we
found no counterexample even against the most optimistic conjecture,
according to which a global solution to (\ref{eq:k}), (\ref{eq:data}),
(\ref{eq:boundary}) exists assuming only that $m$ is a nonnegative
continuous function, and the initial data belong to the ``energy
space'' $V_{0}$.

In this paper we make a first step in the direction of
counterexamples.  We focus on local solutions, and we prove that
Hirosawa's results~\cite{hirosawa-main} are sharp, both in the
strictly hyperbolic and in the weakly hyperbolic case.

In Theorem~\ref{thm:counterex-ndg} and Theorem~\ref{thm:counterex-dg}
we construct indeed local solutions $u(x,t)$ of (\ref{eq:k}),
(\ref{eq:data}), (\ref{eq:boundary}) with very regular initial data,
but such that for every $t>0$ we have that $(u(x,t),u_{t}(x,t))$
belongs to the phase space $V_{\beta}$ if and only if $\beta\leq 1/2$.
In a few words these solutions exhibit an instantaneous
\emph{derivative loss} up to $V_{1/2}$.  In these examples the maximal
regularity admitted for the initial data depends on the continuity
modulus of $m$.  Just to give some examples, in the strictly
hyperbolic case we have the following three situations (see
Example~\ref{ex:ndg}).
\begin{itemize}
	\item If $m$ is continuous we can have derivative loss for quasi
	analytic initial data.  This proves in particular that in
	Nishihara's result~\cite{nishihara} the Lipschitz continuity
	assumption on the nonlinear term cannot be relaxed to continuity,
	even when looking for local solutions.

	\item If $m$ is \aholder\ continuous for some $\alpha\in(0,1)$,
	then we can have derivative loss for initial data in any Gevrey
	space $\G_{s}(\Omega)$ with $s>(1-\alpha)^{-1}$.

	\item If $m$ is \aholder\ continuous for every $\alpha\in(0,1)$,
	then we can have derivative loss for initial data in
	$C^{\infty}(\Omega)$.  This proves in particular that in the
	results stated in (A) (i.e., well posedness in $V_{\beta}$ for
	\emph{every} $\beta\geq 1/2$) the Lipschitz continuity assumption
	on $m$ cannot be relaxed to H\"{o}lder continuity.
\end{itemize}

In the weakly hyperbolic case we have for example the following two
situations (see Example~\ref{ex:dg}).
\begin{itemize}
	\item  If $m$ is \aholder\ continuous for some $\alpha\in(0,1)$,
	then we can have derivative loss for initial data in any Gevrey
	space $\G_{s}(\Omega)$ with $s>1+\alpha/2$.

	\item If $m$ is Lipschitz continuous we can have derivative loss
	for initial data in any Gevrey space $\G_{s}(\Omega)$ with
	$s>3/2$. This proves in particular that in the
	results stated in (A) the strict hyperbolicity cannot be
	relaxed to weak hyperbolicity.
\end{itemize}

Our examples don't exclude that a local solution can always exist in
$V_{1/2}$.  This remains an open problem.

Our approach is based on the theory developed by \textsc{F.\
Colombini}, \textsc{E.\ De Giorgi}, and \textsc{S.\
Spagnolo}~\cite{dgcs}, and then extended by \textsc{F.\ Colombini} and
\textsc{S.\ Spagnolo}~\cite{cs}, and by \textsc{F.\ Colombini},
\textsc{E.\ Jannelli}, and \textsc{S.\ Spagnolo}~\cite{cjs}.  They
considered linear equations with a time dependent
coefficient such as
\begin{equation}
	u_{tt}(t,x)-c(t)\Delta u(t,x)=0.
	\label{eq:linear}
\end{equation}

If $c(t)$ is not Lipschitz continuous or vanishes for $t=0$, they
showed how to construct ``solutions'' of (\ref{eq:linear}) which are
very regular at time $t=0$, but very irregular (not even
distributions) for $t>0$.  

The construction of our counterexamples is divided into two steps.  In
the first step we consider the linear equation (\ref{eq:linear}), and
we modify the parameters in the construction described in
\cite{dgcs,cs,cjs} in order to obtain solutions of (\ref{eq:linear})
which are very regular at time $t=0$, but with a prescribed minimal
regularity (they belong to $V_{\beta}$ if and only if $\beta\leq 1/2$)
for every $t>0$.  These counterexamples for the linear equation,
stated in Proposition~\ref{prop:counterex-ndg} and
Proposition~\ref{prop:counterex-dg}, are maybe interesting in
themselves because they extend the theory developed in
\cite{dgcs,cs,cjs} to coefficients $c(t)$ with arbitrary continuity
modulus.  The assumptions on initial data in these counterexamples are
sharp because they are complementary to those required in the
existence results (both for linear and for Kirchhoff equations).

In the second step we show that these solutions, up to modifying one
dominant Fourier component, are solutions of (\ref{eq:k}) for a
suitable choice of the nonlinearity $m$.

This paper is organized as follows.  In section~\ref{sec:statements}
we rephrase (\ref{eq:k}), (\ref{eq:data}), (\ref{eq:boundary}) as an
abstract evolution problem in a Hilbert space, and we restate
Hirosawa's local existence results in the abstract setting.  In
section~\ref{sec:counterex} we state our counterexamples, which we
prove in section~\ref{sec:proofs}.  For the convenience of the reader
in appendix~\ref{appendix} we sketch the main point of the proof of
the local existence results in the abstract setting.

\setcounter{equation}{0} 
\section{Preliminaries and local existence
results}\label{sec:statements}

\paragraph{Continuity moduli}

Let $\omega:[0,+\infty)\to[0,+\infty)$.  The following assumptions are
recalled throughout the paper.
\begin{enumerate}
	\renewcommand{\labelenumi}{($\omega$\arabic{enumi})}
	\item We have that $\omega\in C^{0}([0,+\infty))$ is an increasing
	function such that $\omega(0)=0$, and
	$\omega(a+b)\leq\omega(a)+\omega(b)$ for every $a\geq 0$ and
	$b\geq 0$.

	\item The function $\sigma\to\sigma/\omega(\sigma)$ is
	nondecreasing.

\end{enumerate}

A function $f:X\to\re$ (where $X\subseteq\re$) is said to be \ocont\
if there exists a constant $L\in\re$ such that
\begin{equation}
	|f(a)-f(b)|\leq
	L\,\omega(|a-b|)
	\hspace{3em}
	\forall a\in X,\ \forall b\in X.
	\label{hp:ocont}
\end{equation}

Two simple properties of continuity moduli are stated in
Lemma~\ref{lemma:omega} of the appendix.  In particular from
(\ref{th:omega-lambda}) it follows that the composition of a Lipschitz
continuous function and an \ocont\ function (in any order) is again an
\ocont\ function.

It is not difficult to see that the set of \ocont\ functions only
depends on the values of $\omega$ in a right neighborhood of 0.  For
this reason, with a little abuse of notation, we often consider
continuity moduli which are defined, and satisfy ($\omega$1) and/or
($\omega$2), in a right neighborhood of 0.  For the same reason when
needed we assume also that $\omega$ is invertible.

\paragraph{Abstract setting for Kirchhoff equations}

Let $H$ be a separable real Hilbert space.  For every $x$ and $y$ in
$H$, let $|x|$ denote the norm of $x$, and let $\langle x,y\rangle$
denote the scalar product of $x$ and $y$.  Let $A$ be an unbounded
linear operator on $H$ with dense domain $D(A)$.  We always assume
that $A$ is self-adjoint and nonnegative, so that the power
$A^{\beta}$ is defined for every $\beta\geq 0$ in a suitable domain
$D(A^{\beta})$.

We consider the second order evolution problem
\begin{equation}
	u''(t)+\m{u(t)}Au(t)=0, 
	\hspace{2em}\foralll t\geq 0,
	\label{pbm:h-eq}
\end{equation}
with initial data
\begin{equation}
	u(0)=u_0,\hspace{3em}u'(0)=u_1.
	\label{pbm:h-data}
\end{equation}
 
It is well known that (\ref{pbm:h-eq}), (\ref{pbm:h-data}) is just an
abstract setting of (\ref{eq:k}), (\ref{eq:data}),
(\ref{eq:boundary}), corresponding to the case where
$H=L^{2}(\Omega)$, and $Au=-\Delta u$, defined for every $u$ in a
suitable domain $D(A)$ depending on the boundary conditions (see
\cite{ap}).

\paragraph{Functional spaces}

For the sake of simplicity, let us assume that $H$ admits a countable
complete orthonormal system $\{e_{k}\}_{k\geq 1}$ made by eigenvectors
of $A$.  We denote the corresponding eigenvalues by $\lambda_{k}^{2}$,
so that $Ae_{k}=\lambda_{k}^{2}e_{k}$ for every $k\geq 1$.

Under this assumption (which in the concrete case corresponds to
bounded domains) we can work with Fourier series.  However, any
definition or statement of this section can be easily extended to the
general setting just by using the spectral theorem for self-adjoint
operators (\cite[Theorem~VIII.4, p.\ 260]{reed}). 

By means of the orthonormal system every $u\in H$ can be written in a
unique way in the form $u=\sum_{k=1}^{\infty}u_{k}e_{k}$, where
$u_{k}=\langle u,e_{k}\rangle$ are the Fourier components of $u$.
Moreover for every $\beta> 0$ we have that
$$A^{\beta}u:=\sum_{k=1}^{\infty}\lambda_{k}^{2\beta}u_{k}e_{k},$$
so that we can consider the quantity
$$\|u\|^{2}_{D(A^{\beta})}:=\sum_{k=1}^{\infty}
\lambda_{k}^{4\beta}u_{k}^{2},$$
and characterize the spaces $D(A^{\beta})$ and $D(A^{\infty})$ as follows
$$D(A^{\beta}):=\left\{u\in
H:\|u\|^{2}_{D(A^{\beta})}<+\infty\right\},
\hspace{3em}
\dainfty:=\bigcap_{\beta>0}^{}D(A^{\beta}).$$

With these notations the phase spaces $V_{\beta}$ are defined as
$V_{\beta}:=D(A^{(\beta+1)/2})\times D(A^{\beta/2})$.

Let now $\varphi:[0,+\infty)\to[1,+\infty)$ be any function.  Then for
every $r>0$ and $\beta>0$ we can set 
\begin{equation}
	\trebar{u}^{2}_{\varphi,r,\beta}:=\sum_{k=1}^{\infty}\lk^{4\beta}
	u_{k}^{2} \exp\left(\strut r\varphi(\lambda_{k})\right),
	\label{defn:trebar}
\end{equation}
and then define the spaces
$$\G_{\varphi,r,\beta}(A):=\left\{u\in
H:\trebar{u}^{2}_{\varphi,r,\beta}<+\infty \right\},
\hspace{3em}
\G_{\varphi,\beta}(A):=\bigcup_{r>0}^{}\G_{\varphi,r,\beta}(A).$$

If two continuous functions $\varphi_{1}(\sigma)$ and
$\varphi_{2}(\sigma)$ coincide for every large enough $\sigma$, then
$\G_{\varphi_{1},r,\beta}(A)=\G_{\varphi_{2},r,\beta}(A)$.  For this
reason, with a little abuse of notation, we consider these spaces
even if $\varphi(\sigma)$ is defined, continuous, and greater than 1
only for large values of $\sigma$.

\paragraph{Hirosawa's results}

We are now ready to recall the main results proved in
\cite{hirosawa-main}, restated in the abstract setting.  Since the
abstract statements do not seem to follow trivially from the ones in
\cite{hirosawa-main}, we sketch the proofs in appendix~\ref{appendix}.

The first result concerns the strictly hyperbolic case (see
Theorem~2.2 in~\cite{hirosawa-main}).

\begin{thmbibl}[Strictly hyperbolic case]\label{hirosawa-ndg}
	Let $\omega:[0,+\infty)\to[0,+\infty)$ be a function satisfying
	($\omega$1).  Let $m:[0,+\infty)\to(0,+\infty)$ be an \ocont\
	function satisfying the strict hyperbolicity condition
	$m(\sigma)\geq\nu>0$ for every $\sigma\geq 0$.
	
	Let $\varphi:[0,+\infty)\to[1,+\infty)$ be a function such that
	\begin{equation}
		\limsup_{\sigma\to +\infty}\frac{\sigma}{\varphi(\sigma)}
		\omega\left(\frac{1}{\sigma}\right)<+\infty.
		\label{hp:phi-ndg-h}
	\end{equation}
	
	Let us assume that
	\begin{equation}
		(u_{0},u_{1})\in\G_{\varphi,r_{0},3/4}(A)\times
		\G_{\varphi,r_{0},1/4}(A)
		\label{hp:data}
	\end{equation}
	for some $r_{0}>0$.
	
	Then there exist $T>0$, and $R>0$ with $RT<r_{0}$ such that
	problem (\ref{pbm:h-eq}), (\ref{pbm:h-data}) admits at least one
	local solution
	\begin{equation}
			u\in C^{1}\left([0,T];\G_{\varphi,r_{0}-Rt,1/4}(A)\right)\cap
			C^{0}\left([0,T];\G_{\varphi,r_{0}-Rt,3/4}(A)\right).
		\label{th:reg-u}
	\end{equation}
\end{thmbibl}

\begin{rmk}
	\begin{em}
		Condition (\ref{th:reg-u}), with the range space depending on
		time, simply means that 
		$$u\in C^{1}\left([0,\tau];\G_{\varphi,r_{0}-R\tau,1/4}(A)\right)\cap
		C^{0}\left([0,\tau];\G_{\varphi,r_{0}-R\tau,3/4}(A)\right)$$
		for every $\tau\in(0,T]$.  
		
		This amounts to say that scales of Hilbert spaces are the
		natural setting for this problem.
	\end{em}
\end{rmk}

\begin{ex}
	\begin{em}
		Let us give some examples in order to clarify the interplay between
		assumptions (\ref{hp:phi-ndg-h}) and (\ref{hp:data}).
		\begin{itemize}
			\item If $\omega(\sigma)=o(1)$ as $\sigma\to 0^{+}$
			(which simply means that $m$ is continuous), then
			(\ref{hp:phi-ndg-h}) holds true with
			$\varphi(\sigma)=\sigma$.  In this case (\ref{hp:data})
			means that $u_{0}$ and $u_{1}$ are analytic, and one has
			re-obtained the classical local existence result for
			analytic data in the strictly hyperbolic case.  In
			(\ref{hp:phi-ndg-h}) one can also take
			$\varphi(\sigma)=\sigma\omega(1/\sigma)$, thus obtaining
			local existence for a larger class of initial data.
		
			\item If $\omega(\sigma)=\sigma^{\alpha}$ for some
			$\alpha\in(0,1)$ (which means that $m$ is H\"{o}lder
			continuous), then (\ref{hp:phi-ndg-h}) holds true with
			$\varphi(\sigma)=\sigma^{1-\alpha}$.  In this case
			(\ref{hp:data}) means that one can take initial data in
			the Gevrey space $\G_{s}(A)$ with $s=(1-\alpha)^{-1}$.  We
			recall that $\G_{s}(A)$ is the space
			$\G_{\varphi,\beta}(A)$ corresponding to
			$\varphi(\sigma)=\sigma^{1/s}$ ($\beta$ is insignificant
			in this case).
		
			\item If $\omega(\sigma)=\sigma|\log\sigma|$ (which means
			that $m$ is log-Lipschitz continuous), then
			(\ref{hp:phi-ndg-h}) holds true with
			$\varphi(\sigma)=\log\sigma$.  In this case
			$\G_{\varphi,r,\beta}(A)=D(A^{\beta+r/4})$.  One can
			therefore take initial data in $D(A^{\infty})$ and obtain
			a solution in the same space, or even initial data in
			$V_{\gamma}$ for some $\gamma>1/2$ and obtain a solution
			with a possible progressive derivative loss (due to the
			term $r_{0}-Rt$ in (\ref{th:reg-u})).
		
			\item If $\omega(\sigma)=\sigma$ (which means that $m$ is
			Lipschitz continuous), then (\ref{hp:phi-ndg-h}) holds
			true with $\varphi(\sigma)\equiv 1$.  Therefore
			(\ref{hp:data}) means that one can take initial data in
			$V_{1/2}$.  This is of course the classical existence
			result in Sobolev spaces.
		\end{itemize}
	\end{em}
\end{ex}

Since we are interested in local solutions, Theorem~\ref{hirosawa-ndg}
can be applied also in the \emph{mildly degenerate} case, namely when
$m$ may vanish but $\m{u_{0}}>0$.  

The following result (see Theorem~2.1 of \cite{hirosawa-main}) concerns
the weakly hyperbolic case, and it is essential to deal with the
\emph{really degenerate} case, i.e., when $\m{u_{0}}=0$.

\begin{thmbibl}[Weakly hyperbolic case]\label{hirosawa-dg}
	Let $\omega:[0,+\infty)\to[0,+\infty)$ be a function satisfying
	($\omega$1).  Let $m:[0,+\infty)\to[0,+\infty)$ be an \ocont\
	function.  
	
	Let $\varphi:[0,+\infty)\to[1,+\infty)$ be a function such that
	\begin{equation}
		\limsup_{\sigma\to +\infty}\sigma
		\left[\varphi\left(\frac{\sigma}{
		\sqrt{\omega(1/\sigma)}}\right)\right]^{-1}
		<+\infty.
		\label{hp:phi-dg-h}
	\end{equation}
	
	Let $(u_{0},u_{1})$ and $r_{0}>0$ be such that (\ref{hp:data})
	holds true.
	
	Then there exist $T>0$, and $R>0$ with $RT<r_{0}$ such that
	problem (\ref{pbm:h-eq}), (\ref{pbm:h-data}) admits at least one
	local solution $u$ satisfying (\ref{th:reg-u}).
\end{thmbibl}

\begin{ex}
	\begin{em}
		Let us give some examples.
		\begin{itemize}
			\item If $\omega(\sigma)=o(1)$ as $\sigma\to 0^{+}$,
			then (\ref{hp:phi-dg-h}) holds true with
			$\varphi(\sigma)=\sigma$.  Once again (\ref{hp:data})
			means that $u_{0}$ and $u_{1}$ are analytic, and one has
			re-obtained the classical local existence result for
			analytic data in the weakly hyperbolic case.
		
			\item If $\omega(\sigma)=\sigma^{\alpha}$ for some
			$\alpha\in(0,1]$, then (\ref{hp:phi-dg-h}) holds true with
			$\varphi(\sigma)=\sigma^{2/(\alpha+2)}$.  Therefore
			(\ref{hp:data}) means that one can take initial data in
			the Gevrey space $\G_{s}(A)$ with $s=1+\alpha/2$.  This is
			true in particular for $\alpha=1$ (i.e., when $m$ is
			Lipschitz continuous).  In this case we have local
			existence for initial data in $\G_{3/2}(A)$.
		\end{itemize}
	\end{em}
\end{ex}

\begin{rmk}
	\begin{em}
		One cannot expect local solutions to be unique if $m$ is not
		Lipschitz continuous (see some examples in \cite{as}).
		However the existence results rely on an a~priori estimate
		(see Proposition~\ref{prop:apriori-est}) which implies in
		particular that under assumptions (\ref{hp:phi-ndg-h}) or
		(\ref{hp:phi-dg-h}) \emph{any} local solution satisfying the
		minimal regularity requirement (\ref{hp:reg-u}) also satisfies
		the stronger condition (\ref{th:reg-u}).  This is a sort of
		propagation of regularity: if the space in (\ref{hp:data}) is
		strictly contained in $D(A^{3/4})\times D(A^{1/4})$, then also
		solutions lie in a scale of spaces which is strictly contained
		in $D(A^{3/4})\times D(A^{1/4})$.  We are going to see that
		this is no more true when condition (\ref{hp:phi-ndg-h}) or
		(\ref{hp:phi-dg-h}) are not satisfied.
	\end{em}
\end{rmk}

\setcounter{equation}{0}
\section{Statements of counterexamples}\label{sec:counterex}

The first counterexample shows the optimality of
Theorem~\ref{hirosawa-ndg} in the non-Lipschitz case.

\begin{thm}[Strictly hyperbolic case]\label{thm:counterex-ndg}
	Let $A$ be a self adjoint linear operator on a Hilbert space $H$.
	Let us assume that there exist a countable (not necessarily
	complete) orthonormal system $\{e_{k}\}_{k\geq 1}$ in $H$, and an
	increasing unbounded sequence $\{\lambda_{k}\}_{k\geq 1}$ of
	positive real numbers such that $Ae_{k}=\lambda_{k}^{2}e_{k}$ for
	every $k\geq 1$.
	
	Let $\omega:[0,+\infty)\to[0,+\infty)$ be a function satisfying
	($\omega$1) and ($\omega$2).
	
	Let $\varphi:[0,+\infty)\to[1,+\infty)$ be a function such that
	\begin{equation}
		\lim_{\sigma\to +\infty}\frac{\sigma}{\varphi(\sigma)}
		\omega\left(\frac{1}{\sigma}\right)=+\infty.
		\label{hp:phi-ndg}
	\end{equation}
	
	Then there exist a function
	$m:[0,+\infty)\to[1/2,3/2]$, a real number $T_{0}>0$, and a
	function $u:[0,T_{0}]\to H$ such that
	\begin{enumerate}
		\renewcommand{\labelenumi}{(\roman{enumi})} 
		\item $m$ is \ocont;
		
		\item $(u(0),u'(0))\in
		\G_{\varphi,r,3/4}(A)\times\G_{\varphi,r,1/4}(A)$ for every
		$r>0$;
	
		\item $u\in C^{1}([0,T_{0}];D(A^{1/4})) \cap
		C^{0}([0,T_{0}];D(A^{3/4}))$ is a solution of
		(\ref{pbm:h-eq});
		
		\item for every $t\in(0,T_{0}]$ we have that
		$(u(t),u'(t))\not\in V_{\beta}$ for every $\beta> 1/2$.
	\end{enumerate}
\end{thm}

\begin{ex}\label{ex:ndg}
	\begin{em}
		Let us consider some examples.
		\begin{itemize}
			\item  The assumptions of Theorem~\ref{thm:counterex-ndg} are
			satisfied if we take
			$$\omega(\sigma)=\frac{1}{|\log\sigma|^{1/2}},
			\hspace{4em}
			\varphi(\sigma)=\frac{\sigma}{\log\sigma}.$$
			
			In this case $m$ is continuous, the initial data are quasi
			analytic, and the solution $u$ has an instantaneous infinite
			derivative loss.
		
			\item Let $\alpha\in(0,1)$.  The assumptions of
			Theorem~\ref{thm:counterex-ndg} are satisfied if we take
			$$\omega(\sigma)=\sigma^{\alpha}, 
			\hspace{4em}
			\varphi(\sigma)=\frac{\sigma^{1-\alpha}}{\log\sigma}.$$
			
			In this case $m$ is $\alpha$-H\"{o}lder continuous, the
			initial data are in the Gevrey space $\G_{s}(A)$ for every
			$s>(1-\alpha)^{-1}$, and the solution $u$ has an
			instantaneous infinite derivative loss.
		
			\item  The assumptions of Theorem~\ref{thm:counterex-ndg} are
			satisfied if we take
			$$\omega(\sigma)=\sigma|\log\sigma|^{3},
			\hspace{4em}
			\varphi(\sigma)=\log^{2}\sigma.$$
			
			In this case $m$ is $\alpha$-H\"{o}lder continuous for
			every $\alpha\in(0,1)$ (but not log-Lipschitz continuous),
			the initial data are in $D(A^{\infty})$, and once again
			the solution $u$ has an instantaneous infinite derivative
			loss.
		\end{itemize}
	\end{em}
\end{ex}

\begin{rmk}\label{rmk:log-lip}
	\begin{em}
		Theorem~\ref{thm:counterex-ndg}, as it is stated, is void in
		the log-Lipschitz case. Indeed when
		$\omega(\sigma)=\sigma|\log\sigma|$ and $\varphi$ satisfies
		(\ref{hp:phi-ndg}), then all initial data satisfying (ii)
		belong to $V_{1/2}$ but not to $V_{\beta}$ for $\beta>1/2$.
		In a certain sense they have no derivatives to loose!
		
		On the other hand, a careful inspection of the proof reveals
		that Theorem~\ref{thm:counterex-ndg} can be improved as
		follows.  Given any function $\psi:[0,+\infty)\to[1,+\infty)$
		such that $\psi(\sigma)\to +\infty$ as $\sigma\to +\infty$ we
		can find a solution satisfying (i), (ii), (iii), and
		\begin{itemize}
			\item[(iv')] for every $t\in(0,T_{0}]$ we have that
			$(u(t),u'(t))\not\in
			\G_{\psi,r,3/4}(A)\times\G_{\psi,r,1/4}(A)$ for every
			$r>0$.
		\end{itemize}
		
		Thus for example we can take 
		$$\omega(\sigma)=\sigma|\log\sigma|,
		\hspace{2em}
		\varphi(\sigma)=\log|\log\sigma|,
		\hspace{2em}
		\psi(\sigma)=\log|\log|\log\sigma||,$$
		and obtain a solution with initial data in
		$\G_{\varphi,r,3/4}(A)\times\G_{\varphi,r,1/4}(A)$ which for
		every positive time doesn't belong even to the weaker scale
		$\G_{\psi,r,3/4}(A)\times\G_{\psi,r,1/4}(A)$.  So also in the
		log-Lipschitz case we have a (infinitesimally small)
		derivative loss.
		
		We spare the reader from this generalization.
	\end{em}
\end{rmk}

The second counterexample shows the optimality of
Theorem~\ref{hirosawa-dg}.

\begin{thm}[Weakly hyperbolic case]\label{thm:counterex-dg}
	Let $A$, $\{e_{k}\}$, $\{\lambda_{k}\}$ be as in
	Theorem~\ref{thm:counterex-ndg}.  
	
	Let $\omega:[0,+\infty)\to[0,+\infty)$ be a function satisfying
	($\omega$1) and ($\omega$2).

	Let $\varphi:[0,+\infty)\to[1,+\infty)$ be a function such that
	\begin{equation}
		\lim_{\sigma\to +\infty}\sigma
		\left[\varphi\left(\frac{\sigma}{
		\sqrt{\omega(1/\sigma)}}\right)\right]^{-1}=+\infty.
		\label{hp:phi-dg}
	\end{equation}
	
	Then there exist a function
	$m:[0,+\infty)\to[0,3/2]$, a real number $T_{0}>0$, and a
	function $u:[0,T_{0}]\to H$ such that
	\begin{enumerate}
		\renewcommand{\labelenumi}{(\roman{enumi})} 
		\item $m$ is \ocont;
		
		\item $(u(0),u'(0))\in
		\G_{\varphi,r,3/4}(A)\times\G_{\varphi,r,1/4}(A)$ for every
		$r>0$;
	
		\item $u\in C^{1}([0,T_{0}];D(A^{1/4})) \cap
		C^{0}([0,T_{0}];D(A^{3/4}))$ is a solution of
		(\ref{pbm:h-eq});
		
		\item for every $t\in(0,T_{0}]$ we have that
		$(u(t),u'(t))\not\in V_{\beta}$ for every $\beta\geq 1$;
		
		\item there exists a sequence $\tau_{k}\to 0^{+}$ such that
		$|A^{(\beta+1)/2}u(\tau_{k})|$ is unbounded for every 
		$\beta> 1/2$.
	\end{enumerate}
\end{thm}

\begin{ex}\label{ex:dg}
	\begin{em}
		Let $\alpha\in(0,1]$.  The assumptions of
		Theorem~\ref{thm:counterex-ndg} are satisfied if we take
		$$\omega(\sigma)=\sigma^{\alpha}, \hspace{4em}
		\varphi(\sigma)=\frac{\sigma^{2/(\alpha+2)}}{\log\sigma}.$$
			
		In this case $m$ is $\alpha$-H\"{o}lder continuous (Lipschitz
		continuous if $\alpha=1$), the initial data are in the Gevrey
		space $\G_{s}(A)$ for every $s>1+\alpha/2$, and the solution
		$u$ has an instantaneous infinite derivative loss.
		
		In particular, in contrast to the strictly hyperbolic case
		(see\ Remark~\ref{rmk:log-lip}),
		Theorem~\ref{thm:counterex-dg} provides examples of infinite
		derivative loss even in the Lipschitz (or log-Lipschitz) case.		
	\end{em}
\end{ex}

The counterexamples of Theorem~\ref{thm:counterex-ndg} and
Theorem~\ref{thm:counterex-dg} originate from two counterexamples for
the linear equation
\begin{equation}
	v''+c(t)Av=0.
	\label{eq:h-lin}
\end{equation}

This equation is well studied in mathematical literature.  It is well
known for example that, if the coefficient $c(t)$ is \ocont, then the
Cauchy problem is well posed in $\G_{\varphi,\beta}(A)$ provided that
$\varphi$ and $\omega$ satisfy (\ref{hp:phi-ndg-h}) in the strictly
hyperbolic case, and (\ref{hp:phi-dg-h}) in the weakly hyperbolic
case.  The main argument is that the approximated energy estimates
introduced in~\cite{dgcs} can be extended word-by-word to arbitrary
continuity moduli as we do in the appendix below.

This result is sharp.  Indeed if $\varphi$ and $\omega$ do not satisfy
(\ref{hp:phi-ndg-h}) or (\ref{hp:phi-dg-h}), hence if they satisfy
(\ref{hp:phi-ndg}) or (\ref{hp:phi-dg}) (see also
Remark~\ref{rmk:liminf} below), then the Cauchy problem is not well
posed in $\G_{\varphi,\beta}(A)$.  In literature we found a lot of
counterexamples for special choices of $\omega$ and $\varphi$ (see for
example~\cite{colombini} and the references quoted therein), but we
didn't find the general counterexamples under assumptions
(\ref{hp:phi-ndg}) or (\ref{hp:phi-dg}).

In the following two propositions we state them in the form which is
needed in the proof of Theorem~\ref{thm:counterex-ndg} and
Theorem~\ref{thm:counterex-dg}.

\begin{prop}[Strictly hyperbolic case]\label{prop:counterex-ndg}
	Let $A$, $\{e_{k}\}$, $\{\lambda_{k}\}$, $\omega$, $\varphi$ be as
	in Theorem~\ref{thm:counterex-ndg}.
	
	Then there exist $c:[0,1]\to[1/2,3/2]$, and $v:[0,1]\to H$ such
	that
	\begin{list}{(SH-\arabic{enumi})}{\usecounter{enumi}\leftmargin
	4em \labelwidth 4em}
		\item $c$ is \ocont;
		
		\item $(v(0),v'(0))\in
		\G_{\varphi,r,3/4}(A)\times\G_{\varphi,r,1/4}(A)$ for every
		$r>0$;
	
		\item \label{th:vcont-sh}
		$v\in C^{1}([0,1];D(A^{1/4})) \cap C^{0}([0,1];D(A^{3/4}))$ is
		a solution of (\ref{eq:h-lin});
		
		\item for every $t\in(0,1]$ we have that
		$(v(t),v'(t))\not\in V_{\beta}$ for every $\beta> 1/2$.
	\end{list}
\end{prop}

\begin{prop}[Weakly hyperbolic case]\label{prop:counterex-dg}
	Let $A$, $\{e_{k}\}$, $\{\lambda_{k}\}$, $\omega$, $\varphi$ be as
	in Theorem~\ref{thm:counterex-dg}.
	
	Then there exist $c:[0,1]\to[0,3/2]$, and $v:[0,1]\to H$ such
	that
	\begin{list}{(WH-\arabic{enumi})}{\usecounter{enumi}\leftmargin
	4em \labelwidth 4em}
		\item $c$ is \ocont;
		
		\item $(v(0),v'(0))\in
		\G_{\varphi,r,3/4}(A)\times\G_{\varphi,r,1/4}(A)$ for every
		$r>0$;
	
		\item $v\in C^{1}([0,1];D(A^{1/4})) \cap
		C^{0}([0,1];D(A^{3/4}))$ is a solution of
		(\ref{eq:h-lin});
		
		\item for every $t\in(0,1]$ we have that
		$(v(t),v'(t))\not\in V_{\beta}$ for every $\beta\geq 1$;
		
		\item there exists a sequence $\tau_{k}\to 0^{+}$ such that
		$|A^{(\beta+1)/2}v(\tau_{k})|$ is unbounded for every
		$\beta> 1/2$.
	\end{list}
\end{prop}

\begin{rmk}
	\begin{em}
		Proposition~\ref{prop:counterex-ndg} and
		Proposition~\ref{prop:counterex-dg} are strongly based on the
		counterexamples shown in \cite{dgcs,cs,cjs}.  On the other
		hand, besides the fact that we deal with arbitrary continuity
		moduli, and arbitrary sequences of eigenvalues, there are some
		differences we would like to point out.
		\begin{itemize}
			\item In \cite{dgcs,cs,cjs} the derivative loss is bigger
			because solutions instantaneously lie outside the space of
			distributions.  Here we need to be more careful since we
			want solutions to lie in $V_{1/2}$ and nothing more.
		
			\item ``Derivative loss'' has a slightly different meaning
			in \cite{dgcs,cs,cjs} and in this paper.  Losing the
			$m$-th derivative in \cite{dgcs,cs,cjs} means that there
			exists a sequence $\tau_{k}\to 0^{+}$ such that the norm
			of $(v(\tau_{k}),v'(\tau_{k}))$ in $V_{m}$ tends to
			$+\infty$.  This of course may happen also if
			$v(\tau_{k})$ and $v'(\tau_{k})$ are in $D(A^{\infty})$
			for every $k$. In other words, what is actually lost is a
			uniform bound on the $m$-th derivative.
			
			In statements (SH-4) and (WH-4), and in the corresponding
			statements of Theorem~\ref{thm:counterex-ndg} and
			Theorem~\ref{thm:counterex-dg}, we lose the $m$-th
			derivative in a \emph{stronger} sense, namely
			$(v(t),v'(t))\not\in V_{m}$ for every $t>0$.  On the
			contrary in statement~(WH-5), and in the corresponding
			statement of Theorem~\ref{thm:counterex-dg}, we are forced
			to lose the last derivatives only in the weaker sense.
			
		\end{itemize}
	\end{em}
\end{rmk}

\begin{rmk}\label{rmk:liminf}
	\begin{em}
		A careful inspection of the proofs shows that the same
		conclusions hold true also if the limit in (\ref{hp:phi-ndg})
		and (\ref{hp:phi-dg}) is replaced by the corresponding limsup 
		computed along the sequence $\{\lambda_{k}\}$.
	\end{em}
\end{rmk}

\setcounter{equation}{0}
\section{Proofs}\label{sec:proofs}

The proof is organized as follows.  In section~\ref{sec:construction}
we construct functions $c(t)$ and $v(t)$ depending on several
parameters.  In Proposition~\ref{prop:c-ocont},
Proposition~\ref{prop:sufficient} and Corollary~\ref{cor:sufficient}
we relate the regularity properties of $c(t)$ and $v(t)$ to suitable
conditions on the parameters.  Then in section~\ref{sec:ndg} we choose
the parameters in order to prove Proposition~\ref{prop:counterex-ndg}.
In section~\ref{sec:dg} we do the same for
Proposition~\ref{prop:counterex-dg}.  Finally in section~\ref{sec:k}
we prove Theorem~\ref{thm:counterex-ndg} and
Theorem~\ref{thm:counterex-dg}.

\subsection{General Construction}\label{sec:construction}

\paragraph{\emph{\textmd{Ingredients}}}

The starting point of the construction are the functions
\begin{eqnarray}
	b(\ep,t) & := & 1-4\ep\sin(2t)-\ep^{2}(1-\cos(2t))^{2},
	\label{defn:bet}  \\
	w(\ep,t) & := & \sin t\cdot\exp\left(\ep\left(t-
	\textstyle{\frac{1}{2}}\sin(2t)\right)\right),
	\label{defn:wet}
\end{eqnarray}
introduced in section~7 of \cite{dgcs}.  

Then the construction is based on seven sequences $\{t_{k}\}$,
$\{s_{k}\}$, $\{\tau_{k}\}$, $\{\eta_{k}\}$, $\{\delta_{k}\}$,
$\{\ep_{k}\}$, $\{a_{k}\}$ satisfying the following assumptions.
\begin{list}{}{\leftmargin 5em \labelwidth 5em}
	\item[(Hp-$t_{k}$)] The sequence $\{t_{k}\}$ is a decreasing
	sequence of positive real numbers (which we think as times) such
	that $t_{0}=1$, and $t_{k}\to 0^{+}$ as $k\to+\infty$.

	\item[(Hp-$s_{k}$)] The sequence $\{s_{k}\}$ is a sequence of
	positive real numbers (which we think as times) such that
	$t_{k}<s_{k}<t_{k-1}$ for every $k\geq 1$.

	\item[(Hp-$\tau_{k}$)] The sequence $\{\tau_{k}\}$ is a sequence of
	positive real numbers (which we think as times) such that
	$t_{k}<\tau_{k}<s_{k}$ for every $k\geq 1$.

	\item[(Hp-$\eta_{k}$)] The sequence $\{\eta_{k}\}$ is an
	increasing subsequence of the sequence $\{\lambda_{k}\}$ of
	the eigenvalues of $A^{1/2}$.  
	
	\item[(Hp-$\delta_{k}$)] The sequence $\{\delta_{k}\}$ is a
	nonincreasing sequence of positive real numbers with
	$\delta_{0}=1$.  Moreover we require that
	$\sqrdk\eta_{k}t_{k}/(2\pi)$ and $\sqrdk\eta_{k}s_{k}/(2\pi)$ are
	integers, and $2\sqrdk\eta_{k}\tau_{k}/\pi$ is an odd integer for
	every $k\geq 1$.

	\item[(Hp-$\ep_{k}$)] The sequence $\{\ep_{k}\}$ is a sequence of
	positive real numbers such that $\ep_{k}\leq 1/16$ for every
	$k\geq 1$.  Moreover we require that $\ep_{k}\delta_{k}\to 0$ and
	that $\{\sqrdk\ep_{k}\eta_{k}\}$ is a nondecreasing sequence.
\end{list} 

No special assumption is required on the sequence $\{a_{k}\}$. We
denote the limit of $\{\delta_{k}\}$ by $\delta_{\infty}$.

\paragraph{\emph{\textmd{Definition and properties of $c(t)$}}}

Let $c:[0,1]\to\re$ be the function defined by $c(0)=\delta_{\infty}$,
and for every $k\geq 1$ 
\begin{equation}
	c(t):=\left\{
	\begin{array}{ll}
		\dk b\left(\ep_{k},\sqrdk\eta_{k}t\right) & \mbox{if }t\in[t_{k},s_{k}],  \\
		\noalign{\vspace{2ex}}
		\displaystyle{\frac{\delta_{k-1}-\delta_{k}}{t_{k-1}-s_{k}}}(t-s_{k})+
		\delta_{k} & \mbox{if }t\in[s_{k},t_{k-1}].
	\end{array}
	\right.
	\label{defn:c}
\end{equation}

Note that in the interval $[s_{k},t_{k-1}]$ the function $c(t)$ is the
affine interpolation between $\dk$ and $\delta_{k-1}$.

From the assumptions on the parameters we have that
$c(t_{k})=c(s_{k})=\delta_{k}$.  Moreover
\begin{equation}
	\dk-8\ep_{k}\delta_{k}\leq c(t)\leq\dk+8\ep_{k}\delta_{k}
	\hspace{3em}
	\forall t\in[t_{k},s_{k}],
	\label{est:c-tksk}
\end{equation}
\begin{equation}
	\dk\leq c(t)\leq \delta_{k-1}
	\hspace{3em}
	\forall t\in[s_{k},t_{k-1}].
	\label{est:c-sktk}
\end{equation}
	
Since $\ep_{k}\dk\to 0$, and $\dk\to\delta_{\infty}$, estimates
(\ref{est:c-tksk}) and (\ref{est:c-sktk}) imply that $c(t)$ is
continuous in $[0,1]$. Since $\ep_{k}\leq 1/16$ we have also that
\begin{equation}
	\frac{1}{2}\delta_{k}\leq c(t)\leq \frac{3}{2}\dk
	\hspace{3em}
	\forall t\in[t_{k},s_{k}],
	\label{est:cepk1-dg}
\end{equation}
and globally
\begin{equation}
	\frac{1}{2}\delta_{\infty}\leq c(t)\leq\frac{3}{2}
	\hspace{3em}
	\forall t\in[0,1].
	\label{est:c-glob}
\end{equation}

Concerning the derivative we have that
\begin{equation}
	|c'(t)|\leq
	16\ep_{k}\eta_{k}\dk^{3/2}
	\hspace{3em}
	\forall t\in(t_{k},s_{k}),
	\label{est:c'1-dg}
\end{equation}
and of course
\begin{equation}
	c'(t)=\frac{\delta_{k-1}-\delta_{k}}{t_{k-1}-s_{k}}
	\hspace{3em}
	\forall t\in(s_{k},t_{k-1}).
	\label{est:c'2-dg}
\end{equation}

The $\omega$-continuity of $c(t)$ in the whole
interval $[0,1]$ can be deduced from the $\omega$-continuity in the
intervals $[t_{k},t_{k-1}]$ provided that the following uniform
estimates hold true (we omit the easy and classical proof).

\begin{lemma}\label{lemma:ocont}
	Let $c:[0,1]\to\re$ be any function.  Let $\{t_{k}\}$ be any
	sequence satisfying (Hp-$t_{k}$).
	
	Then $c$ is \ocont\ in $[0,1]$ if and only if there exists a
	constant $L$ such that
	\begin{enumerate}
		\renewcommand{\labelenumi}{(\roman{enumi})}
		\item $|c(t_{i})-c(t_{j})|\leq L\,\omega(|t_{i}-t_{j}|)$ for
		every pair of positive integers $i$ and $j$;
	
		\item $|c(a)-c(b)|\leq L\,\omega(|a-b|)$ for every $k\geq 1$,
		and every $a$ and $b$ in $[t_{k},t_{k-1}]$.\qed
	\end{enumerate}
\end{lemma}

Applying the lemma we find the following sufficient condition for the
for the $\omega$-continuity of $c(t)$.

\begin{prop}\label{prop:c-ocont}
	The function $c(t)$ defined in (\ref{defn:c}) is \ocont\ in
	$[0,1]$ if
	\begin{equation}
		\sup_{i<j}\frac{\delta_{i}-\delta_{j}}{\omega(t_{i}-t_{j})}+
		\sup_{k\geq
		1}\frac{\delta_{k-1}-\delta_{k}}{\omega(t_{k-1}-s_{k})}+
		\sup_{k\geq 1}
		\frac{\ep_{k}\delta_{k}}{\omega\left(
		2\pi/(\eta_{k}\sqrdk)\right)}
		<+\infty.
		\label{est:c-ocont}
	\end{equation}
\end{prop}

\prf We apply Lemma~\ref{lemma:ocont}.  Assumption (i) is satisfied
because the first supremum in (\ref{est:c-ocont}) is finite.  In order
to verify assumption~(ii) we consider separately the subintervals
$[t_{k},s_{k}]$ and $[s_{k},t_{k-1}]$.  Thanks to ($\omega$2) in
$[s_{k},t_{k-1}]$ we have that $$\frac{|c(a)-c(b)|}{\omega(|a-b|)}=
\frac{\delta_{k-1}-\delta_{k}}{t_{k-1}-s_{k}}\cdot
\frac{|a-b|}{\omega(|a-b|)}\leq
\frac{\delta_{k-1}-\delta_{k}}{\omega(t_{k-1}-s_{k})}.$$

This is bounded uniformly in $k$ because the second supremum in
(\ref{est:c-ocont}) is finite.

In $[t_{k},s_{k}]$ the function $c(t)$ is periodic, hence we can limit
ourselves to consider points $a$ and $b$ with $|a-b|$ less or equal
than a period $2\pi/(\eta_{k}\sqrdk)$.  By (\ref{est:c'1-dg}) and
($\omega$2) we therefore have that
$$\frac{|c(a)-c(b)|}{\omega(|a-b|)}=
\left|\frac{c(a)-c(b)}{a-b}\right|\cdot
\frac{|a-b|}{\omega(|a-b|)}\leq 16\ep_{k}\eta_{k}\delta_{k}^{3/2}\cdot
\frac{2\pi}{\eta_{k}\sqrdk}\cdot \frac{1}{\omega\left(
2\pi/(\eta_{k}\sqrdk)\right)}.$$

This is bounded uniformly in $k$ because the third supremum in
(\ref{est:c-ocont}) is finite.  \qed

\paragraph{\emph{\textmd{Definition and properties of $v_{k}(t)$}}}

Let $v_{k}:[0,1]\to\re$ be the solution of the linear ordinary
differential equation
\begin{equation}
	v_{k}''(t)+\eta_{k}^{2}c(t)v_{k}(t)=0,
	\label{eq:ODE}
\end{equation}
with ``initial'' data
\begin{equation}
	v_{k}(t_{k})=0,
	\hspace{3em}
	v_{k}'(t_{k})=\eta_{k}\sqrdk.
	\label{eq:PDE-data}
\end{equation}

In order to estimate $v_{k}$ we consider the usual ``Kovalevskian'' energy
\begin{equation}
	E_{k}(t) :=  |v_{k}'(t)|^{2}+\eta_{k}^{2}|v_{k}(t)|^{2},
	\label{defn:ek}  \\
\end{equation}
and the usual hyperbolic energy
\begin{equation}
	F_{k}(t)  :=  |v_{k}'(t)|^{2}+\eta_{k}^{2}c(t)|v_{k}(t)|^{2}.
	\label{defn:fk}
\end{equation}

It is simple to show that their time derivatives are given by
\begin{equation}
	E_{k}'(t)=2\eta_{k}^{2}(1-c(t))v_{k}(t)v_{k}'(t);
	\label{eq:deriv-ek}
\end{equation}
\begin{equation}
	F_{k}'(t)=\eta_{k}^{2}c'(t)|v_{k}(t)|^{2}.
	\label{eq:deriv-fk}
\end{equation}

Now we estimate $E_{k}(t)$ in four cases.

\subparagraph{\emph{\textmd{Case 1: $t\in[t_{k},s_{k}]$}}}

In this interval we have the explicit expression
$$v_{k}(t)=\exp\left(-\ep_{k}\eta_{k}\sqrdk t_{k}\right)\cdot 
w\left(\ep_{k},\eta_{k}\sqrdk t\right),$$
which is indeed the main motivation of definitions (\ref{defn:bet})
and (\ref{defn:wet}).  In particular we have the following equalities
\begin{equation}
	E_{k}(t_{k})=F_{k}(t_{k})=\dk\eta_{k}^{2},
	\label{eq:ektk-dg}
\end{equation}
\begin{equation}
	E_{k}(s_{k})=F_{k}(s_{k})=\dk\eta_{k}^{2}
	\exp\left(2\ep_{k}\eta_{k}\sqrdk(s_{k}-t_{k})\right),
	\label{eq:eksk-dg}
\end{equation}
and the estimate
\begin{equation}
	E_{k}(t)\leq 3\eta_{k}^{2}\exp\left(2\ep_{k}\eta_{k}\sqrdk
	(s_{k}-t_{k})\right)
	\hspace{2em}
	\forall t\in[t_{k},s_{k}].
	\label{est:ekik}
\end{equation}

From the explicit expression, and our assumption that
$2\sqrdk\eta_{k}\tau_{k}/\pi$ is an odd integer, we have also that
\begin{equation}
	\left|v_{k}(\tau_{k})\right|=\exp\left(
	\ep_{k}\eta_{k}\sqrt{\dk}(\tau_{k}-t_{k})\right).
	\label{eq:vktauk}
\end{equation}

\subparagraph{\emph{\textmd{Case 2: $t\in[s_{k},t_{k-1}]$}}}

In this interval $c'(t)\geq 0$, and therefore from (\ref{eq:deriv-fk})
we have that 
$$F_{k}'(t)\leq\frac{c'(t)}{c(t)}F_{k}(t) 
\hspace{2em}
\forall t\in[s_{k},t_{k-1}].$$

Integrating in $[s_{k},t]$ we find that
\begin{equation}
	F_{k}(t)\leq F_{k}(s_{k})
	\exp\left(\int_{s_{k}}^{t}\frac{c'(s)}{c(s)}\,ds\right)
	= F_{k}(s_{k})
	\exp\left(\log\frac{c(t)}{c(s_{k})}\right)
	= F_{k}(s_{k})\cdot\frac{c(t)}{\dk}
	\label{est:fkski}
\end{equation}
for every $t\in[s_{k},t_{k-1}]$. Using (\ref{eq:eksk-dg}), and the fact
that $c(t)\leq\delta_{k-1}$ in this interval, we obtain in particular 
that
\begin{equation}
	F_{k}(t)\leq\delta_{k-1}\eta_{k}^{2}\exp\left(2\ep_{k}\eta_{k}\sqrdk
	(s_{k}-t_{k})\right)
	\hspace{2em}
	\forall t\in[s_{k},t_{k-1}].
	\label{est:fksk+dg}
\end{equation}

On the other hand we have that $F_{k}'(t)\geq 0$ in this interval,
hence
\begin{equation}
	F_{k}(t)\geq F_{k}(s_{k})=\dk\eta_{k}^{2}
	\exp\left(2\ep_{k}\eta_{k}\sqrdk(s_{k}-t_{k})\right)
	\hspace{2em}
	\forall t\in[s_{k},t_{k-1}].
	\label{est:fkskd+dg}
\end{equation}

Since $c(t)\leq 3/2$ we have also that
$$E_{k}(t)\leq\frac{3}{2}\cdot\frac{1}{c(t)}F_{k}(t)
\hspace{2em}
\forall t\in(0,1],$$
which combined with (\ref{est:fkski}) and (\ref{eq:eksk-dg}) gives
\begin{equation}
	E_{k}(t)\leq\frac{3}{2}\eta_{k}^{2}\exp\left(
	2\ep_{k}\eta_{k}\sqrdk(s_{k}-t_{k})\right)
	\hspace{2em}
	\forall t\in[s_{k},t_{k-1}].
	\label{est:eksk+dg}
\end{equation}

\subparagraph{\emph{\textmd{Case 3: $t\in[0,t_{k}]$}}}

From (\ref{eq:deriv-ek}) and (\ref{est:c-glob}) we have that
$$\left|E_{k}'(t)\right|\leq\eta_{k}|1-c(t)|E_{k}(t) \leq
\eta_{k}E_{k}(t)
\hspace{3em}
\forall t\in[0,1].$$

Integrating in $[t,t_{k}]$ we obtain that
$$E_{k}(t_{k})\exp(-\eta_{k}(t_{k}-t))\leq E_{k}(t)\leq
E_{k}(t_{k})\exp(\eta_{k}(t_{k}-t))
\hspace{2em}
\forall t\in[0,t_{k}],$$
hence by (\ref{eq:ektk-dg}) we conclude that
\begin{equation}
	\dk\eta_{k}^{2}\exp(-\eta_{k}t_{k})\leq E_{k}(t)\leq
	\dk\eta_{k}^{2}\exp(\eta_{k}t_{k})
	\hspace{2em}
	\forall t\in[0,t_{k}].
	\label{est:ek0-dg}
\end{equation}

\subparagraph{\emph{\textmd{Case 4: $t\in[t_{k-1},1]$}}}

From (\ref{eq:deriv-fk}) we have that
$$F_{k}'(t)=\frac{c'(t)}{c(t)}\cdot\eta_{k}^{2}c(t)|v_{k}(t)|^{2}
\leq\frac{|c'(t)|}{c(t)}F_{k}(t),$$
hence
\begin{equation}
	F_{k}(t)\leq F_{k}(t_{k-1})\exp\left(
	\int_{t_{k-1}}^{1}\frac{|c'(s)|}{c(s)}\,ds\right)
	\hspace{2em}
	\forall t\in[t_{k-1},1].
	\label{est:fk-exp}
\end{equation}

Now we observe that
\begin{equation}
	\int_{t_{k-1}}^{1}\frac{|c'(s)|}{c(s)}\,ds=
	\sum_{i=1}^{k-1}\int_{t_{i}}^{s_{i}}\frac{|c'(s)|}{c(s)}\,ds+
	\sum_{i=1}^{k-1}\int_{s_{i}}^{t_{i-1}}\frac{|c'(s)|}{c(s)}\,ds,
	\label{eq:int-c'c}
\end{equation}
and we estimate the two summands in the right hand side.  From
(\ref{est:c'1-dg}) and (\ref{est:cepk1-dg}) we have that
$$\frac{|c'(s)|}{c(s)}\leq
\frac{16\ep_{i}\eta_{i}\delta_{i}^{3/2}}{\delta_{i}/2}=
32\ep_{i}\eta_{i}\sqrt{\delta_{i}}
\hspace{2em}
\forall s\in[t_{i},s_{i}].$$

By (Hp-$\ep_{k}$) the sequence $\ep_{k}\eta_{k}\sqrdk$ is
nondecreasing.  Therefore we have that
$$\sum_{i=1}^{k-1}\int_{t_{i}}^{s_{i}}\frac{|c'(s)|}{c(s)}\,ds\leq
32\ep_{k-1}\eta_{k-1}\sqrt{\delta_{k-1}}
\sum_{i=1}^{k-1}|s_{i}-t_{i}|\leq
32\ep_{k-1}\eta_{k-1}\sqrt{\delta_{k-1}}.$$

On the other hand, in the interval $[s_{i},t_{i-1}]$ we have that
$c'(s)>0$, hence 
$$\int_{s_{i}}^{t_{i-1}}\frac{|c'(s)|}{c(s)}\,ds=
\int_{s_{i}}^{t_{i-1}}\frac{c'(s)}{c(s)}\,ds=
\log\frac{c(t_{i-1})}{c(s_{i})}=
\log\frac{\delta_{i-1}}{\delta_{i}},$$
and therefore (we recall that $\delta_{0}=1$)
$$\sum_{i=1}^{k-1}\int_{s_{i}}^{t_{i-1}}\frac{|c'(s)|}{c(s)}\,ds=
\sum_{i=1}^{k-1}\log\frac{\delta_{i-1}}{\delta_{i}}=
\log\left(\prod_{i=1}^{k-1}\frac{\delta_{i-1}}{\delta_{i}}\right)=
\log\frac{1}{\delta_{k-1}}.$$

Coming back to (\ref{eq:int-c'c}) and (\ref{est:fk-exp}) we have
proved that
$$F_{k}(t)\leq F_{k}(t_{k-1})\frac{1}{\delta_{k-1}}
\exp\left(32\ep_{k-1}\eta_{k-1}\sqrt{\delta_{k-1}}\right).$$

Taking (\ref{est:fksk+dg}) into account we finally obtain that
\begin{equation}
	F_{k}(t)\leq 
	\eta_{k}^{2}\exp\left(2\ep_{k}\eta_{k}\sqrdk(s_{k}-t_{k})+
	32\ep_{k-1}\eta_{k-1}\sqrt{\delta_{k-1}}\right)
	\hspace{2em}
	\forall t\in[t_{k-1},1].
	\label{est:fktu-dg}
\end{equation}

In an analogous way we can prove that
$$F_{k}(t)\geq F_{k}(t_{k-1})\exp\left(-
\int_{t_{k-1}}^{1}\frac{|c'(s)|}{c(s)}\,ds\right)\geq
F_{k}(t_{k-1})\delta_{k-1}\exp\left(-
32\ep_{k-1}\eta_{k-1}\sqrt{\delta_{k-1}}\right),$$
which by (\ref{est:fkskd+dg}) gives
\begin{equation}
	F_{k}(t)\geq \delta_{k}\delta_{k-1}
	\eta_{k}^{2}\exp\left(2\ep_{k}\eta_{k}\sqrdk(s_{k}-t_{k})-
	32\ep_{k-1}\eta_{k-1}\sqrt{\delta_{k-1}}\right)
	\hspace{2em}
	\forall t\in[t_{k-1},1].
	\label{est:fktd-dg}
\end{equation}

In order to obtain estimates on $E_{k}(t)$ we just remark that
$$\frac{\delta_{k-1}}{2}\leq c(t)\leq \frac{3}{2}	
\hspace{2em}
\forall t\in[t_{k-1},1]$$
(here we used that the sequence $\delta_{k}$ is nonincreasing), so
that
$$\frac{2}{3}F_{k}(t)\leq E_{k}(t)\leq
\frac{3}{\delta_{k-1}}F_{k}(t).$$

Therefore (\ref{est:fktu-dg}) and (\ref{est:fktd-dg}) imply the
following estimates for every $t\in[t_{k-1},1]$
\begin{equation}
	E_{k}(t)\leq 
	\frac{3\eta_{k}^{2}}{\delta_{k-1}}
	\exp\left(2\ep_{k}\eta_{k}\sqrdk(s_{k}-t_{k})+
	32\ep_{k-1}\eta_{k-1}\sqrt{\delta_{k-1}}\right),
	\label{est:ektu-dg}
\end{equation}
\begin{equation}
	E_{k}(t)\geq \frac{2}{3}\delta_{k}\delta_{k-1}
	\eta_{k}^{2}\exp\left(2\ep_{k}\eta_{k}\sqrdk(s_{k}-t_{k})-
	32\ep_{k-1}\eta_{k-1}\sqrt{\delta_{k-1}}\right).
	\label{est:ektd-dg}
\end{equation}

\paragraph{\emph{\textmd{Definition and properties of $v(t)$}}}

Let $e_{i_{k}}$ denote the eigenvector of $A$ corresponding to the
eigenvalue $\eta_{k}^{2}$. We set
\begin{equation}
	v(t):=\sum_{k=1}^{\infty}a_{k}v_{k}(t)e_{i_{k}}.
	\label{defn:v}
\end{equation}

At this level of generality this definition is purely formal, because
we have no information about the convergence of the series.  All
regularity properties of $v(t)$ are related to convergence properties
of suitable series.  In particular for every $t\in[0,1]$, and every
$\beta\geq 0$, $\varphi:[0,+\infty)\to[1,+\infty)$, and $r>0$ we have
that
\begin{eqnarray}
	(v(t),v'(t))\in V_{\beta}&
	\Longleftrightarrow &
	\sum_{k=1}^{\infty}a_{k}^{2}\eta_{k}^{2\beta}E_{k}(t)<+\infty,
	\label{impl:v-sobolev}  \\
	\hspace{-2em}
	(v(t),v'(t))\in\G_{\varphi,r,3/4}(A)\times\G_{\varphi,r,1/4}(A) &
	\Longleftrightarrow &
	\sum_{k=1}^{\infty}a_{k}^{2}\eta_{k}E_{k}(t)
	\exp\left(\strut r\varphi(\eta_{k})\right)<+\infty.
	\label{impl:v-gevrey}
\end{eqnarray}

Combining these implications with our estimates on $E_{k}(t)$ we
obtain sufficient conditions in terms of the parameters for the
regularity or non regularity of $v(t)$.

\begin{prop}\label{prop:sufficient}
	Let $v(t)$ be the function defined in (\ref{defn:v}).
	\begin{enumerate}
		\renewcommand{\labelenumi}{(\arabic{enumi})}
		\item For every $r>0$ we have that $(v(0),v'(0))
		\in\G_{\varphi,r,3/4}(A)\times\G_{\varphi,r,1/4}(A)$ if
		\begin{equation}
			\sum_{k=1}^{\infty}a_{k}^{2}\dk\eta_{k}^{3}\exp\left(
			\eta_{k}t_{k}+\strut
			r\varphi(\eta_{k})\right)<+\infty.
			\label{est:v0-gevrey-dg}
		\end{equation}
	
		\item For every $t>0$ and $\beta\geq 0$ we have that
		$(v(t),v'(t))\in V_{\beta}$ if
		\begin{equation}
			\sum_{k=1}^{\infty}\frac{a_{k}^{2}\eta_{k}^{2\beta+2}}{\delta_{k-1}}
			\exp\left(2\ep_{k}\eta_{k}\sqrdk(s_{k}-t_{k})+
			32\ep_{k-1}\eta_{k-1}\sqrt{\delta_{k-1}}\right)
			<+\infty.
			\label{est:vt-sobolev-dg}
		\end{equation}
	
		\item For every $t>0$ and $\beta\geq 0$ we have that
		$(v(t),v'(t))\not\in V_{\beta}$ if
		\begin{equation}
			\sum_{k=1}^{\infty}a_{k}^{2}\delta_{k}\delta_{k-1}\eta_{k}^{2\beta+2}
			\exp\left(2\ep_{k}\eta_{k}\sqrdk(s_{k}-t_{k})-
			32\ep_{k-1}\eta_{k-1}\sqrt{\delta_{k-1}}\right)
			=+\infty.
			\label{est:vt-nosobolev-dg}
		\end{equation}
	
		\item We have that $v\in C^{1}([0,1];D(A^{1/4})) \cap
		C^{0}([0,1];D(A^{3/4}))$ if (\ref{est:v0-gevrey-dg}) holds
		true with $r=0$, and (\ref{est:vt-sobolev-dg}) holds true with
		$\beta=1/2$.
	
		\item  For every $\beta\geq 0$ we have that the sequence
		$\left|A^{(\beta+1)/2}v(\tau_{k})\right|$ is unbounded if
		\begin{equation}
			\lim_{k\to +\infty}a_{k}^{2}\eta_{k}^{2\beta+2}\exp\left(
			2\ep_{k}\eta_{k}\sqrdk(\tau_{k}-t_{k})\right)=+\infty.
			\label{est:sup-infty}
		\end{equation}
	\end{enumerate}
\end{prop}

\prf
Statement (1) follows combining (\ref{impl:v-gevrey}) and
(\ref{est:ek0-dg}) with $t=0$.

Let us examine statements~(2) and~(3).  For $t\in(0,1]$ we have that
$t\geq t_{k-1}$ for every large enough $k$, and therefore the relevant
estimates are (\ref{est:ektu-dg}) and (\ref{est:ektd-dg}).  At this
point statement~(2) follows combining (\ref{impl:v-sobolev}) with
(\ref{est:ektu-dg}), statement (3) follows by combining
(\ref{impl:v-sobolev}) with (\ref{est:ektd-dg}).

Let us consider now statement (4).  A sufficient condition in order to
prove that $v(t)$ has the required time regularity is that
\begin{equation}
	\sum_{k=1}^{\infty}a_{k}^{2}\eta_{k}
	\sup_{t\in[0,1]}E_{k}(t)<+\infty.
	\label{est:series}
\end{equation}

Now for every $k\geq 1$ we have that
$$\sup_{t\in[0,1]}E_{k}(t)\leq
\dk\eta_{k}^{2}\exp(\eta_{k}t_{k})+
\frac{3\eta_{k}^{2}}{\delta_{k-1}}
\exp\left(2\ep_{k}\eta_{k}\sqrdk(s_{k}-t_{k})+
32\ep_{k-1}\eta_{k-1}\sqrt{\delta_{k-1}}\right).$$

This estimate indeed easily follows from (\ref{est:ek0-dg}) if $t\leq
t_{k}$, from (\ref{est:ekik}) or (\ref{est:eksk+dg}) if $t_{k}\leq
t\leq t_{k-1}$, from (\ref{est:ektu-dg}) if $t_{k-1}\leq t\leq 1$.

Therefore (\ref{est:series}) follows from (\ref{est:v0-gevrey-dg})
with $r=0$ and (\ref{est:vt-sobolev-dg}) with $\beta=1/2$.

Let us finally prove statement (5). By (\ref{eq:vktauk}) we have that
$$\left|A^{(\beta+1)/2}v(\tau_{k})\right|^{2}\geq
a_{k}^{2}\eta_{k}^{2\beta+2}|v_{k}(\tau_{k})|^{2}=
a_{k}^{2}\eta_{k}^{2\beta+2}\exp\left(
2\ep_{k}\eta_{k}\sqrdk(\tau_{k}-t_{k})\right),$$
so that (\ref{est:sup-infty}) implies that the left hand side is
unbounded.\qed

\paragraph{\emph{\textmd{Choice of $a_{k}$}}}

Let us set
\begin{equation}
	a_{k}^{2}:=\frac{\delta_{k-1}}{k^{2}\eta_{k}^{3}}
	\exp\left(-2\ep_{k}\eta_{k}\sqrdk(s_{k}-t_{k})-
	32\ep_{k-1}\eta_{k-1}\sqrt{\delta_{k-1}}\right).
	\label{defn:ak-dg}
\end{equation}

Conditions (\ref{est:v0-gevrey-dg}) through (\ref{est:sup-infty}) can
thus be restated in terms of the remaining parameters.  We obtain the
following result (proof is trivial).

\begin{cor}\label{cor:sufficient}
	Let $v(t)$ be the function defined in (\ref{defn:v}), with $a_{k}$
	defined by (\ref{defn:ak-dg}).
	\begin{enumerate}
		\renewcommand{\labelenumi}{(\arabic{enumi})} \item For every
		$r>0$ we have that
		$(v(0),v'(0))\in\G_{\varphi,r,3/4}(A)\times\G_{\varphi,r,1/4}(A)$
		if
		\begin{equation}
			\sum_{k=1}^{\infty}\frac{\delta_{k}\delta_{k-1}}{k^{2}}
			\exp\left(
			\eta_{k}t_{k}+r\varphi(\eta_{k})
			-2\ep_{k}\eta_{k}\sqrdk(s_{k}-t_{k})-
			32\ep_{k-1}\eta_{k-1}\sqrt{\delta_{k-1}}\right)<+\infty.
			\label{est:v0-gevrey-red}
		\end{equation}
	
		\item For every $t>0$ and $\beta\geq 0$ we have that
		$(v(t),v'(t))\in V_{\beta}$ if
		\begin{equation}
			\sum_{k=1}^{\infty}
			\frac{1}{k^{2}}\eta_{k}^{2\beta-1}< +\infty.
			\label{est:vt-sobolev-red}
		\end{equation}
	
		\item For every $t>0$ and $\beta\geq 0$ we have that
		$(v(t),v'(t))\not\in V_{\beta}$ if
		\begin{equation}
			\sum_{k=1}^{\infty}
			\frac{\delta_{k}\delta_{k-1}^{2}}{k^{2}}\eta_{k}^{2\beta-1}
			\exp\left(-64\ep_{k-1}\eta_{k-1}\sqrt{\delta_{k-1}}\right)
			=+\infty.
			\label{est:vt-nosobolev-red}
		\end{equation}
	
		\item We have that $v\in C^{1}([0,1];D(A^{1/4})) \cap
		C^{0}([0,1];D(A^{3/4}))$ if (\ref{est:v0-gevrey-red}) holds
		true with $r=0$ (note that (\ref{est:vt-sobolev-red}) is
		trivial for $\beta=1/2$).
	
		\item  For every $\beta\geq 0$ we have that the sequence
		$\left|A^{(\beta+1)/2}v(\tau_{k})\right|$ is unbounded if
		\begin{equation}
			\lim_{k\to
			+\infty}\frac{\delta_{k-1}}{k^{2}}\eta_{k}^{2\beta-1}
			\exp\left(2\ep_{k}\eta_{k}\sqrdk(\tau_{k}-s_{k})-
			32\ep_{k-1}\eta_{k-1}\sqrt{\delta_{k-1}}\right)=+\infty.
			\label{est:sup-infty-red}
		\end{equation}
	\end{enumerate}
\end{cor}

\subsection{Proof of Proposition
\ref{prop:counterex-ndg}}\label{sec:ndg}

We define $c(t)$ and $v(t)$ according to the general construction.  In
order to choose the parameters we first set $t_{0}=\delta_{0}=1$.
Then for every $k\geq 1$ we set
\begin{equation}
	t_{k}:=\frac{2\pi}{\eta_{k}},
	\hspace{2em}
	s_{k}:=\left[\frac{t_{k-1}}{2t_{k}}\right]t_{k},
	\hspace{2em}
	\delta_{k}:=1,
	\hspace{2em}
	\ep_{k}:=\omega\left(\frac{1}{\eta_{k}}\right),
	\label{defn:tsek}
\end{equation}
where the brackets denote the integer part, and the sequence
$\eta_{k}$ is defined as follows:
\begin{itemize}
	\item  $\eta_{1}$ is any element of the sequence $\{\lambda_{k}\}$
	which is greater than $4\pi$ and such that $\omega(1/\eta_{1})\leq
	1/16$; 

	\item for every $k\geq 2$ the term $\eta_{k}$ is recursively
	defined as any element of the sequence $\{\lambda_{k}\}$
	satisfying the following three inequalities
	\begin{eqnarray}
		\eta_{k} & \geq & 4\eta_{k-1},
		\label{eq:eta1}  \\
		\frac{\eta_{k}}{\varphi(\eta_{k})}
		\omega\left(\frac{1}{\eta_{k}}\right) & \geq & \frac{k}{t_{k-1}},
		\label{eq:eta2}  \\
		\eta_{k} & \geq & \exp(k+k\ep_{k-1}\eta_{k-1}).
		\label{eq:eta3}  
	\end{eqnarray}
\end{itemize}

We didn't define $\tau_{k}$ because we don't need it in this proof.

The first thing we have to verify is that $\eta_{k}$ is well defined.
The right hand sides of (\ref{eq:eta1}) through (\ref{eq:eta3}) only
depend on the values of the sequences with index $k-1$.  Since the
sequence $\{\lambda_{k}\}$ is unbounded, we can find $\eta_{k}$ with
the required properties if we show that the left hand sides tend to
$+\infty$ as $\eta_{k}\to +\infty$.  This is trivial in the case of
(\ref{eq:eta1}) and (\ref{eq:eta3}), and it follows from assumption
(\ref{hp:phi-ndg}) in the case of (\ref{eq:eta2}).

The second thing we need is that assumptions (Hp-$t_{k}$),
(Hp-$s_{k}$), (Hp-$\eta_{k}$), (Hp-$\delta_{k}$), and (Hp-$\ep_{k}$)
are satisfied (they are required for the entire construction to work).
By (\ref{eq:eta1}) the sequence $\eta_{k}$ is increasing and tends to
$+\infty$.  These facts are enough to prove all properties, except the
monotonicity of $\ep_{k}\eta_{k}\sqrdk$, which however follows from
($\omega$2).

We are now ready to prove the conclusions of
Proposition~\ref{prop:counterex-ndg}.

\paragraph{\textmd{\emph{Conclusion (SH-1)}}}

First of all we remark that $1/2\leq c(t)\leq 3/2$. This is due to
(\ref{est:c-glob}) and our definition of $\delta_{k}$.

In order to prove the $\omega$-continuity of $c(t)$ we apply
Proposition~\ref{prop:c-ocont}.  The first and second supremum in
(\ref{est:c-ocont}) are trivially zero.  Moreover by the monotonicity
of $\omega$ we have that 
$$\frac{\ep_{k}\delta_{k}}{\omega\left( 2\pi/(\eta_{k}\sqrdk)\right)}
=\frac{\ep_{k}}{\omega(t_{k})}=
\frac{\omega(t_{k}/(2\pi))}{\omega(t_{k})}\leq 1,$$
hence also the third supremum is finite.

\paragraph{\textmd{\emph{Conclusion (SH-2)}}}

Let us examine (\ref{est:v0-gevrey-red}).  The series converges for
every $r>0$ (and also for $r=0$) if we show that the argument of the
exponential function is bounded from above.  Since $\eta_{k}t_{k}$ and
$\ep_{k}\eta_{k}t_{k}$ are bounded, it is enough to prove that
$r\varphi(\eta_{k})\leq 2\ep_{k}\eta_{k}s_{k}$ for every large enough
$k$.  In turn this is true if we show that
\begin{equation}
	\lim_{k\to
	+\infty}\frac{\ep_{k}\eta_{k}s_{k}}{\varphi(\eta_{k})}=+\infty.
	\label{lim-c3}
\end{equation}

To this end we first use the definition of $t_{k}$ and (\ref{eq:eta1})
to obtain that
$$\frac{s_{k}}{t_{k-1}}=
\left[\frac{t_{k-1}}{2t_{k}}\right]\frac{t_{k}}{t_{k-1}}\geq
\frac{1}{2}-\frac{t_{k}}{t_{k-1}}=
\frac{1}{2}-\frac{\eta_{k-1}}{\eta_{k}}\geq
\frac{1}{4}.$$

Then using the definition of $s_{k}$ and $\ep_{k}$, and estimate
(\ref{eq:eta2}), we deduce that
$$\frac{\ep_{k}\eta_{k}s_{k}}{\varphi(\eta_{k})}=
\frac{\eta_{k}}{\varphi(\eta_{k})}\omega\left(\frac{1}{\eta_{k}}\right)
s_{k}\geq k\cdot \frac{s_{k}}{t_{k-1}}\geq \frac{k}{4},$$
which implies (\ref{lim-c3}).

\paragraph{\textmd{\emph{Conclusion (SH-3)}}}

By statement (4) of Corollary~\ref{cor:sufficient} it is enough to
prove (\ref{est:v0-gevrey-red}) with $r=0$.  We did this in the
previous paragraph.

\paragraph{\textmd{\emph{Conclusion (SH-4)}}}

Let us examine (\ref{est:vt-nosobolev-red}).  By
(\ref{eq:eta3}) we have that 
$$\frac{1}{k^{2}}\eta_{k}^{2\beta-1}
\exp(-64\ep_{k-1}\eta_{k-1})\geq
\frac{1}{k^{2}}e^{(2\beta-1)k}\exp\left(\strut
[(2\beta-1)k-64]\ep_{k-1}\eta_{k-1}\right).$$

It is therefore clear that the series cannot converge if $\beta>1/2$.
\qed

\subsection{Proof of Proposition \ref{prop:counterex-dg}}\label{sec:dg}

Let $g:(0,+\infty)\to(0,+\infty)$ be defined by
$g(\sigma):=\sqrt{\sigma}\omega^{-1}(\sigma)$, where $\omega^{-1}$ is
the inverse function of $\omega$.  It is easy to see that $g$ is
invertible.  Its inverse function, which we denote by $h(\sigma)$, is
increasing and vanishes as $\sigma\to 0^{+}$.

Now we define $c(t)$ and $v(t)$ according to the general construction.
In order to choose the parameters we first set $t_{0}=\delta_{0}=1$.
Then for every $k\geq 1$ we set
\begin{equation}
	t_{k}:=\frac{2\pi}{\eta_{k}\sqrdk},
	\hspace{1em}
	s_{k}:=\left[\frac{t_{k-1}}{2t_{k}}\right]t_{k},
	\hspace{1em}
	\tau_{k}:=s_{k}-\frac{t_{k}}{4},
	\hspace{1em}
	\ep_{k}:=\frac{1}{16},
	\hspace{1em}
	\delta_{k}:=h\left(\frac{1}{\eta_{k}}\right),
	\label{defn:tsek-dg}
\end{equation}
where the brackets denote the integer part, and the sequence
$\eta_{k}$ is defined as follows:
\begin{itemize}
	\item $\eta_{1}$ is any element of the sequence $\{\lambda_{k}\}$
	with $h(1/\eta_{1})<1$, and $\eta_{1}\sqrt{h(1/\eta_{1})}>4\pi$;

	\item for every $k\geq 2$ the term $\eta_{k}$ is recursively
	defined as any element of the sequence $\{\lambda_{k}\}$ such that
	$\eta_{k}\geq\eta_{k-1}+1$ and
	\begin{eqnarray}
		\eta_{k}\left[h\left(\frac{1}{\eta_{k}}\right)\right]^{1/2} & \geq & 
		4\eta_{k-1}\sqrt{\delta_{k-1}},
		\label{eq:eta1-dg}  \\
		\eta_{k}h\left(\frac{1}{\eta_{k}}\right)   & \geq & 
		\frac{k}{t_{k-1}},
		\label{eq:eta2-dg}  \\
		\frac{\eta_{k}}{\varphi(\eta_{k})}
		\left[h\left(\frac{1}{\eta_{k}}\right)\right]^{1/2}   & \geq & 
		\frac{k}{t_{k-1}},
		\label{eq:eta3-dg}  \\
		\eta_{k}h\left(\frac{1}{\eta_{k}}\right) & \geq & 
		\frac{k^{2}}{\delta_{k-1}^{2}}\exp\left(
		64\ep_{k-1}\eta_{k-1}\sqrt{\delta_{k-1}}\right),
		\label{eq:eta4-dg}  \\
		\log\eta_{k} & \geq & 
		k\left(2\eta_{k-1}\sqrt{\delta_{k-1}}-\log\delta_{k-1}
		+2\log k\right).
		\label{eq:eta5-dg}
	\end{eqnarray}
\end{itemize}

The first thing we need is that $\eta_{k}$ is well defined.  As in the
proof of Proposition~\ref{prop:counterex-ndg} it is enough to show
that the left hand sides of (\ref{eq:eta1-dg}) through
(\ref{eq:eta5-dg}) tend to $+\infty$ as $\eta_{k}\to +\infty$.  This
is trivial in the case of (\ref{eq:eta5-dg}). As for
(\ref{eq:eta1-dg}), (\ref{eq:eta2-dg}), (\ref{eq:eta4-dg}), with the
variable change $z=1/g(\omega(\sigma))$ we have that
$$\lim_{z\to +\infty}
z\left[h\left(\frac{1}{z}\right)\right]^{\gamma}=
\lim_{\sigma\to 0^{+}}
\frac{\left[\omega(\sigma)\right]^{\gamma}}{g(\omega(\sigma))}=
\lim_{\sigma\to 0^{+}}
\left[\omega(\sigma)\right]^{\gamma-3/2}\frac{\omega(\sigma)}{\sigma}.$$

By ($\omega$2) the function $\omega(\sigma)/\sigma$ is bounded from
below in a right neighborhood of $0$.  It follows that the limit is
$+\infty$ for every $\gamma<3/2$.

Let us finally consider the left hand side of
(\ref{eq:eta3-dg}).  Thanks to assumption (\ref{hp:phi-dg}) and the
variable change $z=1/g(\omega(1/\sigma))$ we have that
\begin{eqnarray*}
	\lim_{z\to +\infty} \frac{z}{\varphi(z)}
	\left[h\left(\frac{1}{z}\right)\right]^{1/2} & = &
	\lim_{\sigma\to +\infty} \frac{1}{g(\omega(1/\sigma))}\cdot
	\frac{1}{\varphi(1/g(\omega(1/\sigma)))}\cdot
	\left[\omega(1/\sigma)\right]^{1/2} \\
	 & = & \lim_{\sigma\to +\infty} \sigma
	 \left[\varphi\left(\frac{\sigma}{
	 \sqrt{\omega(1/\sigma)}}\right)\right]^{-1}=+\infty.
\end{eqnarray*}

The second thing we need is that assumptions (Hp-$t_{k}$),
(Hp-$s_{k}$), (Hp-$\tau_{k}$), (Hp-$\eta_{k}$), (Hp-$\delta_{k}$), and
(Hp-$\ep_{k}$) are satisfied.  Assumption (\ref{eq:eta1-dg}) is
equivalent to say that $\eta_{k}\sqrdk\geq
4\eta_{k-1}\sqrt{\delta_{k-1}}$, hence also $t_{k}\leq t_{k-1}/4$.
All the other properties are now almost trivial.

We are now ready to prove the conclusions of
Proposition~\ref{prop:counterex-dg}.

\paragraph{\textmd{\emph{Conclusion (WH-1)}}}

First of all we remark that $0\leq c(t)\leq 3/2$ (with $c(t)=0$ if and
only if $t=0$).  This is due to (\ref{est:cepk1-dg}),
(\ref{est:c-glob}), and our definition of $\delta_{k}$.

In order to prove the $\omega$-continuity of $c(t)$ we apply
Proposition~\ref{prop:c-ocont}.  Let us consider the first supremum in
(\ref{est:c-ocont}). By definition of $t_{k}$ and $\delta_{k}$ we have
that $\delta_{k}=\omega(t_{k}/(2\pi))$. Since the function $\omega$ is
itself \ocont, for $j>i$ we have that
$$\delta_{i}-\delta_{j}=
\omega\left(\frac{t_{i}}{2\pi}\right)-
\omega\left(\frac{t_{j}}{2\pi}\right)\leq
\omega\left(\frac{t_{i}-t_{j}}{2\pi}\right)\leq
\omega(t_{i}-t_{j}),$$
which implies that the first supremum in (\ref{est:c-ocont}) is finite.

Since clearly  $s_{k}/t_{k-1}\leq 1/2$, we have that
$t_{k-1}-s_{k}\geq t_{k-1}/2$, hence
$$\delta_{k-1}-\delta_{k}\leq
\delta_{k-1}=\omega\left(\frac{t_{k-1}}{2\pi}\right)\leq
\omega\left(\frac{t_{k-1}}{2}\right)\leq
\omega\left(t_{k-1}-s_{k}\right),$$
which implies that the second supremum in (\ref{est:c-ocont}) is
finite.  Finally 
$$\delta_{k}=\omega\left(\frac{t_{k}}{2\pi}\right)\leq
\omega\left(t_{k}\right)=
\omega\left(\frac{2\pi}{\eta_{k}\sqrdk}\right),$$
hence also the third supremum is trivially finite.

\paragraph{\textmd{\emph{Conclusion (WH-2)}}}

Let us examine (\ref{est:v0-gevrey-red}).  The series converges for
every $r>0$ if we show that the argument of the exponential function
is bounded from above. Since $\ep_{k}\eta_{k}\sqrdk t_{k}$ is
constant, it is enough to prove that
$$\eta_{k}t_{k}+r\varphi(\eta_{k})\leq
2\ep_{k}\eta_{k}\sqrdk s_{k}$$
whenever $k$ is large enough.  Since $\ep_{k}$ is constant this is
true if we show that
\begin{equation}
	\lim_{k\to +\infty}
	\frac{\sqrdk s_{k}}{t_{k}}=+\infty,
	\hspace{3em}
	\lim_{k\to +\infty}
	\frac{\eta_{k}\sqrdk s_{k}}{\varphi(\eta_{k})}=+\infty.
	\label{lim-c3-dg}
\end{equation}

Since we already observed that $t_{k}/t_{k-1}\leq 1/4$, we have also
that 
\begin{equation}
	\frac{s_{k}}{t_{k-1}}=
	\left[\frac{t_{k-1}}{2t_{k}}\right]\frac{t_{k}}{t_{k-1}}\geq
	\frac{1}{2}-\frac{t_{k}}{t_{k-1}}\geq
	\frac{1}{4}.
	\label{eq:sktk-dg}
\end{equation}

Combining with (\ref{eq:eta2-dg}) we obtain that
$$\frac{\sqrdk s_{k}}{t_{k}}=
\frac{\dk\eta_{k}s_{k}}{2\pi}=
\frac{1}{2\pi}\cdot h\left(\frac{1}{\eta_{k}}\right)\eta_{k}\cdot
\frac{s_{k}}{t_{k-1}}\cdot t_{k-1}\geq\frac{k}{8\pi},$$
which implies the first limit in (\ref{lim-c3-dg}).  Using
(\ref{eq:sktk-dg}) and (\ref{eq:eta3-dg}) we have that
$$\frac{\eta_{k}\sqrdk s_{k}}{\varphi(\eta_{k})}=
\frac{\eta_{k}}{\varphi(\eta_{k})}
\left[h\left(\frac{1}{\eta_{k}}\right)\right]^{1/2}
\frac{s_{k}}{t_{k-1}}\cdot t_{k-1}\geq\frac{k}{4},$$
which implies the second limit in (\ref{lim-c3-dg}).

\paragraph{\textmd{\emph{Conclusion (WH-3)}}}

As in the strictly hyperbolic case it follows from the previous
paragraph.

\paragraph{\textmd{\emph{Conclusion (WH-4)}}}

Let us examine (\ref{est:vt-nosobolev-red}).  Assumption
(\ref{eq:eta4-dg}) implies that the general term of the series is
greater or equal than $\eta_{k}^{2\beta-2}$.  Therefore the series
diverges for every $\beta\geq 1$.

\paragraph{\textmd{\emph{Conclusion (WH-5)}}}

Let us finally examine (\ref{est:sup-infty-red}).  By
(\ref{eq:eta5-dg}) each term in the sequence is greater or equal than
$$\exp\left([(2\beta-1)k-1]
\left(2\eta_{k-1}\sqrt{\delta_{k-1}}-\log\delta_{k-1} +2\log
k\right)-\frac{\pi}{16}\right),$$
hence its limit is $+\infty$ whenever $\beta>1/2$.\qed

\subsection{Proof of Theorem \ref{thm:counterex-ndg} and Theorem
\ref{thm:counterex-dg}}\label{sec:k}

For shortness's sake we only prove Theorem~\ref{thm:counterex-ndg}
(the proof of Theorem~\ref{thm:counterex-dg} is completely analogous).

Let $c(t)$ and $v(t)$ be as in Proposition~\ref{prop:counterex-ndg}.
Since properties of $v(t)$ don't depend on a single Fourier component,
we can always assume that the component of $v(t)$ with respect to
$e_{1}$ is identically zero.

Now we consider the solution $w:[0,1]\to\re$ of the Cauchy problem
$$w''(t)+\lambda_{1}^{2}c(t)w(t)=0, 
\hspace{3em} 
w(0)=1,\ w'(0)=1,$$
and for every $\ep>0$ we set
$$\uep(t):=w(t)e_{1}+\ep v(t)
\hspace{2em}
\forall t\in[0,1].$$
	
The function $\uep:[0,1]\to H$ is again a solution of
(\ref{eq:h-lin}), and it has the same regularity and non regularity
properties of $v(t)$.

Let $\psi_{\ep}(t):=\auq{\uep(t)}$.  Since the component of $v(t)$ with
respect to $e_{1}$ is zero we have that 
$$\psi_{\ep}'(t)=2\langle
A^{3/4}\uep(t),A^{1/4}\uep'(t)\rangle= 2\lambda_{1}^{2} w(t)w'(t)+2\ep
\langle A^{3/4}v(t),A^{1/4}v'(t)\rangle.$$

From (SH-3) we have that $\langle A^{3/4}v(t),A^{1/4}v'(t)\rangle$ is
bounded in $[0,1]$.  Since $w(0)w'(0)=1$ there exists $T_{0}\in(0,1]$
and $\ep_{0}>0$ such that $$\psi_{\ep_{0}}'(t)\geq\frac{1}{2}
\hspace{2em}
\forall t\in[0,T_{0}].$$

It follows that
$\psi_{\ep_{0}}:[0,T_{0}]\to[\psi_{\ep_{0}}(0),\psi_{\ep_{0}}(T_{0})]$
is invertible, and therefore we can define $m:[0,+\infty)\to[1/2,3/2]$
by 
$$m(\sigma):=\left\{
\begin{array}{ll}
	c(0) & \mbox{if } \sigma\leq\psi_{\ep_{0}}(0), \\
	\noalign{\vspace{0.5ex}}
	c\left(\psi_{\ep_{0}}^{-1}(\sigma)\right) & \mbox{if } 
	\psi_{\ep_{0}}(0)\leq\sigma\leq\psi_{\ep_{0}}(T_{0}),   \\
	\noalign{\vspace{0.5ex}}
	c(T_{0}) & \mbox{if } \sigma\geq\psi_{\ep_{0}}(T_{0}).
\end{array}\right.$$

The function $m$ is \ocont\ because $\psi_{\ep_{0}}^{-1}$ is Lipschitz
continuous. Moreover
$$c(t)=\m{u_{\ep_{0}}(t)}
\hspace{2em}
\forall t\in[0,T_{0}].$$

Therefore for this choice of the nonlinearity $m$ the function
$u(t):=u_{\ep_{0}}(t)$ is a solution of (\ref{pbm:h-eq}) in the
interval $[0,T_{0}]$, and it exhibits the required derivative
loss.\qed
\appendix
\setcounter{equation}{0}
\section{Proof of Theorem~\ref{hirosawa-ndg} and Theorem~\ref{hirosawa-dg}
(sketch)}\label{appendix}

It is well known that local existence for Kirchhoff equations can be
easily proved in different ways (for example Galerkin approximations)
provided that solutions are known to satisfy an a priori estimate
yielding the uniform continuity of the nonlinear term.  In this
appendix we prove this estimate under the assumptions of
Theorem~\ref{hirosawa-ndg} and Theorem~\ref{hirosawa-dg}.  The
statement is the following.

\begin{prop}\label{prop:apriori-est}
	Let $\omega$, $m$, $\varphi$, $u_{0}$, $u_{1}$, $r_{0}$ be as in
	Theorem~\ref{hirosawa-ndg} or Theorem~\ref{hirosawa-dg}.
	
	Then there exist positive real numbers $T$, $H$, $R$, with
	$RT<r_{0}$, such that every solution
	\begin{equation}
		u\in C^{0}\left([0,T];D(A^{3/4})\right)\cap
			C^{1}\left([0,T];D(A^{1/4})\right)
		\label{hp:reg-u}
	\end{equation}
	of problem (\ref{pbm:h-eq}), (\ref{pbm:h-data}) satisfies
	\begin{equation}
		|A^{1/4}u'(t)|^{2}+|A^{3/4}u(t)|^{2}\leq H
		\quad\quad\forall t\in[0,T],
	\label{th:apriori-est}
	\end{equation}
	and
	\begin{equation}
		u\in C^{1}\left([0,T];\G_{\varphi,r_{0}-Rt,1/4}(A)\right)\cap
		C^{0}\left([0,T];\G_{\varphi,r_{0}-Rt,3/4}(A)\right).
		\label{th:apriori-reg}
	\end{equation}
\end{prop}

The constants $T$, $H$, $R$ depend only on $\omega$, $m$, and on the
norms of $u_{0}$ and $u_{1}$ in the spaces $\G_{\varphi,r_{0},3/4}(A)$
and $\G_{\varphi,r_{0},1/4}(A)$, respectively.

In the proof of Proposition~\ref{prop:apriori-est} we need two simple 
properties of continuity moduli and convolutions (see Lemma~4.1 and
Lemma~4.2 in \cite{gg:k-manfrinosawa}).

\begin{lemma}\label{lemma:omega}
	Let $\omega:[0,+\infty)\to[0,+\infty)$ be a continuity modulus
	satisfying ($\omega$1).
	
	Then we have that
	\begin{eqnarray}
		 & \omega(\lambda \sigma)\leq(1+\lambda)\omega(\sigma)
		 \quad\quad\forall\lambda\geq 0,\ \forall \sigma\geq 0; & 
		\label{th:omega-lambda}  \\
		\noalign{\vspace{1ex}}
		 & \displaystyle{1+\frac{1}{\omega(\sigma)}
		 \leq\left(1+\frac{1}{\omega(1)}\right)
		 \left(1+\frac{1}{\sigma}\right)}
		 \quad\quad\forall \sigma>0.& 
		\label{th:omega-3}
	\end{eqnarray}
\end{lemma}

\begin{lemma}\label{lemma:conv}
	Let $\rho:\re\to [0,+\infty)$ be a function of class $C^{\infty}$,
	with support contained in $[-1,1]$, and integral equal to 1.
	
	Let $a>0$, and let $f:[0,a]\to\re$ be a continuous function.  Let
	us extend $f(x)$ to the whole real line by setting $f(x)=f(0)$ for
	every $x\leq 0$, and $f(x)=f(a)$ for every $x\geq a$.  
	
	For every $\ep>0$ let us set 
	$$f_{\ep}(x):=\int_{\re}^{}f(x+\ep
	s)\rho(s)\,ds \quad\quad
	\forall x\in\re.$$
	
	Then $f_{\ep}\in C^{\infty}(\re)$, and $f_{\ep}$ has the following
	properties.
	\begin{enumerate}
		\renewcommand{\labelenumi}{(\arabic{enumi})}
		\item If $\mu_{1}\leq f(x)\leq\mu_{2}$ for every $x\in[0,a]$,
		then $\mu_{1}\leq f_{\ep}(x)\leq\mu_{2}$ for every
		$x\in\re$ and every $\ep>0$.

		\item Let $\omega$ be a continuity modulus. Let us assume that 
		\begin{equation}
			|f(x)-f(y)|\leq H_{0}\,\omega(|x-y|) \quad\quad \forall x\in[0,a],\
			\forall y\in[0,a],
			\label{hp:o-cont}
		\end{equation}
		for some $H_{0}\geq 0$.  Then there exists a constant
		$\gamma_{0}$ (depending on $\rho$, but independent on $\ep$,
		$H_{0}$, $f$) such that
		$$|f_{\ep}(x)-f(x)|\leq\gamma_{0}H_{0}\,\omega(\ep) \quad\quad
			\forall x\in\re,\ \forall\ep>0,$$
		$$|f_{\ep}'(x)|\leq\gamma_{0}H_{0}\,\displaystyle{\frac{\omega(\ep)}{\ep}}
		\quad\quad
			\forall x\in\re,\ \forall\ep>0.$$

	\end{enumerate}
\end{lemma}

\paragraph{A priori estimates in the strictly hyperbolic case}

Let us introduce some constants. From the assumptions on $m$ and
$\varphi$ we know that there exist $L$ and $\Lambda$ such that
\begin{equation}
	|m(\sigma_{1})-m(\sigma_{2})|\leq L\,\omega(|\sigma_{1}-\sigma_{2}|)
	\quad\quad\forall\sigma_{1}\geq 0,\ \forall\sigma_{2}\geq 0,
	\label{defn:L}
\end{equation}
\begin{equation}
	\sigma\omega(1/\sigma)\leq\Lambda\varphi(\sigma)
	\quad\quad\forall\sigma> 0.
	\label{defn:Lambda}
\end{equation}

Let $\gamma_{0}$ be the constant appearing in Lemma~\ref{lemma:conv},
and let
$$\gamma_{1}:=\left(\m{u_{0}}+1\right)\cdot\max\left\{1,\nu^{-1}\right\},$$
$$H:=\gamma_{1}\left(\trebar{u_{1}}^{2}_{\varphi,r_{0},1/4}+
\trebar{u_{0}}^{2}_{\varphi,r_{0},3/4}\right)+1,$$
$$R:=\gamma_{0}L\Lambda(H+1)\left(
\frac{1}{\nu}+\frac{1}{\sqrt{\nu}}\right).$$

Let $T>0$ be such that
\begin{equation}
	\omega(T)\leq\frac{1}{L(H+1)},
	\hspace{4em}
	T\leq\frac{r_{0}}{R}.
	\label{defn:T-ndg}
\end{equation}

We claim that (\ref{th:apriori-est}) and (\ref{th:apriori-reg}) hold
true for these values of $T$, $H$, and $R$.  To this end we introduce
\begin{equation}
	S:=\sup\left\{\tau\leq T:|A^{1/4}u'(t)|^{2}+|A^{3/4}u(t)|^{2}\leq H
	\;\;\forall t\in[0,\tau]\right\}.
	\label{defn:S}
\end{equation}

We have that $S>0$ because $\gamma_{1}\geq 1$ and 
$$\|u_{1}\|^{2}_{D(A^{1/4})}+\|u_{0}\|^{2}_{D(A^{3/4})}\leq
\trebar{u_{1}}^{2}_{\varphi,r_{0},1/4}+
\trebar{u_{0}}^{2}_{\varphi,r_{0},3/4}.$$

In order to prove (\ref{th:apriori-est}) we only need to show that
$S=T$.  So let us assume by contradiction that $S<T$.  By the
maximality of $S$ we have that necessarily
\begin{equation}
	|A^{1/4}u'(S)|^{2}+|A^{3/4}u(S)|^{2}= H.
	\label{eq:S-nec}
\end{equation}

Moreover for every $t\in[0,S]$ we have that
\begin{equation}
	\left|\frac{\mathrm{d}}{\mathrm{d}t}|A^{1/2}u(t)|^{2}\right|=
	2\left|\langle A^{3/4}u(t),A^{1/4}u'(t)
	\rangle\right|\leq H.
	\label{eq:c'}
\end{equation}

Therefore if we set 
\begin{equation}
	c(t):=\m{u(t)},
	\label{defn:ct}
\end{equation}
by (\ref{defn:L}), (\ref{eq:c'}), and (\ref{th:omega-lambda}) we have that
$$\left|c(t)-c(s)\right|\leq L\,\omega\left(\left|
|A^{1/2}u(t)|^{2}-|A^{1/2}u(s)|^{2}\right|\right)\leq
L(H+1)\,\omega(|t-s|)$$
for every $t$ and $s$ in $[0,S]$.  In particular from the strict
hyperbolicity and the first inequality in (\ref{defn:T-ndg}) we obtain
that
$$\nu\leq c(t)\leq c(0)+L(H+1)\omega(t)\leq\m{u_{0}}+1
\quad\quad\forall t\in[0,S].$$

Let us extend $c(t)$ outside the
interval $[0,S]$ as in Lemma~\ref{lemma:conv}, and let us set
\begin{equation}
	\cep(t):=\int_{\re}^{}c(t+\ep s)\rho(s)\,ds \quad\quad \forall
	t\in\re.
	\label{defn:cep}
\end{equation}

Since estimate (\ref{hp:o-cont}) holds true with $H_{0}:=L(H+1)$, from
statements (1) and (2) of Lemma~\ref{lemma:conv} we deduce that
\begin{equation}
	\nu\leq\cep(t)\leq\m{u_{0}}+1
	\quad\quad
	\forall t\in\re,\ \forall\ep>0,  
	\label{est:cep-nu-mu}
\end{equation}
\begin{equation}
	 |\cep(t)-c(t)|\leq\gamma_{0}L(H+1)\,\omega(\ep)
	\quad\quad
	\forall t\in\re,\ \forall\ep>0,  
	\label{est:cep-c}
\end{equation}
\begin{equation}
	|\cep'(t)|\leq\gamma_{0}L(H+1)\,
	\displaystyle{\frac{\omega(\ep)}{\ep}}
	\quad\quad
	\forall t\in\re,\ \forall\ep>0. 
	\label{est:cep'}
\end{equation}

Let us consider the Fourier components $u_{k}(t)$ of $u(t)$, and let
us set 
\begin{equation}
	E_{k,\ep}(t):=|u_{k}'(t)|^{2}+\lk^{2}\cep(t)|u_{k}(t)|^{2}.
	\label{defn:ekep}
\end{equation}

An easy computation shows that
\begin{eqnarray*}
	E_{k,\ep}'(t) & = & \cep'(t)\lk^{2}|u_{k}(t)|^{2}+
	2\lk^{2}(\cep(t)-c(t))u_{k}(t)u_{k}'(t)\\
	 & \leq & \frac{|\cep'(t)|}{\cep(t)}E_{k,\ep}(t)+
	 \lk\frac{|\cep(t)-c(t)|}{\sqrt{\cep(t)}}E_{k,\ep}(t),
\end{eqnarray*}
hence by (\ref{est:cep-nu-mu}), (\ref{est:cep-c}), and
(\ref{est:cep'}) we obtain that
\begin{equation}
	E_{k,\ep}'(t)\leq \gamma_{0}L(H+1)\left(
	\frac{1}{\nu}\frac{\omega(\ep)}{\ep}+
	\frac{1}{\sqrt{\nu}}\lk\omega(\ep)\right)E_{k,\ep}(t)
	\quad\quad\forall t\in[0,S].
	\label{est:ekep}
\end{equation}

Let us consider now the eigenvalues $\lk>0$, and let us set
$\ep_{k}:=1/\lk$.  By (\ref{defn:Lambda}) we have that
$$\frac{\omega(\ep_{k})}{\ep_{k}}=\lk\omega(\ep_{k})=
\lk\omega\left(\frac{1}{\lk}\right)\leq\Lambda\varphi(\lk),$$
hence by (\ref{est:ekep})
$$E_{k,\ep_{k}}'(t)\leq \gamma_{0}L(H+1)\left(
\frac{1}{\nu}+\frac{1}{\sqrt{\nu}}\right)
\Lambda\varphi(\lk) E_{k,\ep_{k}}(t)=
R\varphi(\lk)E_{k,\ep_{k}}(t).$$

Integrating this differential inequality we find that
$$E_{k,\ep_{k}}(t)\leq E_{k,\ep_{k}}(0)
\exp\left(R\varphi(\lk)t\right)
\quad\quad\forall t\in[0,S].$$

Combining this inequality with (\ref{est:cep-nu-mu}) we obtain that
\begin{eqnarray}
	|u_{k}'(t)|^{2}+\lk^{2}|u_{k}(t)|^{2} & \leq & 
	\max\left\{1,\nu^{-1}\right\}E_{k,\ep_{k}}(t)
	\nonumber  \\
	 & \leq & \max\left\{1,\nu^{-1}\right\}\cdot
	 (c(0)+1)
	 \left(|u_{1k}|^{2}+\lk^{2} |u_{0k}|^{2}\right)
	 \exp\left(Rt\varphi(\lk)\right)
	 \nonumber  \\
	 & = & \gamma_{1}\left(|u_{1k}|^{2}+\lk^{2}|u_{0k}|^{2}\right)
	 \exp\left(Rt\varphi(\lk)\right),
	 \label{est:uk-data}
\end{eqnarray}
where $u_{0k}$ and $u_{1k}$ denote the Fourier components of $u_{0}$
and $u_{1}$, respectively.

Recalling that $RT\leq r_{0}$, and the definition of $H$, we 
have that
\begin{eqnarray*}
	|A^{1/4}u'(t)|^{2}+|A^{3/4}u(t)|^{2} & =
	& \sum_{k}\lk\left( |u_{k}'(t)|^{2}+\lk^{2}|u_{k}(t)|^{2}\right)
	\\
	 & \leq & \gamma_{1} \sum_{k}\lk
	 \left(|u_{1k}|^{2}+\lk^{2}|u_{0k}|^{2}\right)
	 \exp\left(r_{0}\varphi(\lk)\right) \\
	 & = &
	 \gamma_{1}\left(\trebar{u_{1}}^{2}_{\varphi,r_{0},1/4}+
	 \trebar{u_{0}}^{2}_{\varphi,r_{0},3/4}\right)<H,
\end{eqnarray*}
for every $t\in[0,S]$.  This inequality with $t=S$ contradicts
(\ref{eq:S-nec}) and thus proves (\ref{th:apriori-est}).

A sufficient condition in order to prove (\ref{th:apriori-reg}) is
that 
$$\mathcal{E}(\tau):=\sum_{k}\lk\exp\left((r_{0}-R\tau)\varphi(\lk)\right)
\sup_{t\in[0,\tau]}\left( |u_{k}'(t)|^{2}+\lk^{2}|u_{k}(t)|^{2}\right)
<+\infty$$
for every $\tau\in(0,T]$.  On the other hand from (\ref{est:uk-data})
we have that
\begin{eqnarray*}
	\mathcal{E}(\tau) & \leq & \gamma_{1} \sum_{k}\lk
	 \left(|u_{1k}|^{2}+\lk^{2}|u_{0k}|^{2}\right)
	 \exp\left(r_{0}\varphi(\lk)\right) \\
	 & = & \gamma_{1}\left(\trebar{u_{1}}^{2}_{\varphi,r_{0},1/4}+
	 \trebar{u_{0}}^{2}_{\varphi,r_{0},3/4}\right).
\end{eqnarray*}

%

\paragraph{A priori estimates in the weakly hyperbolic case}

Let $L$ and $\gamma_{0}$ be the constants appearing in (\ref{defn:L})
and in Lemma~\ref{lemma:conv}.  Let $\Lambda$ be such that
\begin{equation}
	\sigma\leq\Lambda\varphi\left(\frac{\sigma}{
		\sqrt{\omega(1/\sigma)}}\right)
		\quad\quad
		\forall\sigma> 0.
	\label{defn:Lambda-dg}
\end{equation}

Let us set
\begin{eqnarray*}
	 & \displaystyle{\gamma_{2}:=1+\frac{1}{\omega(1)}}, \hspace{3em}
	 \gamma_{3}:=\omega(1)+ \max\left\{\omega(\sigma):0\leq
	 \sigma\sqrt{\omega(\sigma)}\leq 1\right\}, &  \\
	 \noalign{\vspace{0.5ex}}
	 & \displaystyle{\gamma_{4}:=2\gamma_{2}}
	 \left(\gamma_{3}+\m{u_{0}}+1\right), &  \\
	 \noalign{\vspace{0.5ex}}
	 &
	H:=\gamma_{4}\displaystyle{\left(1+\frac{2\Lambda}{r_{0}}\right)}
	\left(\trebar{u_{1}}^{2}_{\varphi,r_{0},1/4}+
	\trebar{u_{0}}^{2}_{\varphi,r_{0},3/4}\right)+1, &  \\
	\noalign{\vspace{0.5ex}} & \gamma_{5}:=\gamma_{0}L(H+1)+1,
	\hspace{2em}
	\gamma_{6}:=2\gamma_{5}\max\left\{\Lambda,1+\sqrt{\omega(1)}\right\},
	\hspace{2em}R:=2\gamma_{6}.  &
\end{eqnarray*}

Let $T>0$ be such that
\begin{equation}
	\omega(T)\leq\frac{1}{L(H+1)},
	\hspace{3em}
	T\leq\frac{r_{0}}{2\gamma_{6}}.
	\label{defn:T-dg}
\end{equation}

We claim that in the weakly hyperbolic case (\ref{th:apriori-est}) and
(\ref{th:apriori-reg}) hold true with these new values of $T$, $H$,
and $R$.  To this end we define once again $S$ according to
(\ref{defn:S}), and we assume by contradiction that $S<T$, so that
(\ref{eq:S-nec}) holds true.

Let us consider the function $c(t)$ defined according to
(\ref{defn:ct}), let us extend it outside the interval $[0,S]$ as in
Lemma~\ref{lemma:conv}, and let us set
\begin{equation}
	\cep(t):=\omega(\ep)+\int_{\re}^{}c(t+\ep s)\rho(s)\,ds
	\quad\quad\forall t\in\re.
	\label{defn:cep-dg}
\end{equation}

Arguing as in the strictly hyperbolic case we find that
$$\left|c(t)-c(s)\right| \leq L(H+1)\omega(|t-s|)$$ 
for every $t$ and $s$ in $[0,S]$. This estimate, together with the
first inequality in (\ref{defn:T-dg}), gives that
\begin{equation}
	c(t)\leq c(0)+L(H+1)\omega(t)\leq\m{u_{0}}+1
	\quad\quad\forall t\in[0,S].
	\label{est:c-dg}
\end{equation}

Moreover from statement (2) of Lemma~\ref{lemma:conv} we deduce that
\begin{equation}
	 |\cep(t)-c(t)|\leq\left(1+\gamma_{0}L(H+1)\right)\omega(\ep)=
	 \gamma_{5}\,\omega(\ep),  
	\label{est:cep-c-w}
\end{equation}
\begin{equation}
	|\cep'(t)|\leq\gamma_{0}L(H+1)
	\displaystyle{\frac{\omega(\ep)}{\ep}}\leq
	\gamma_{5}\,
	\displaystyle{\frac{\omega(\ep)}{\ep}}.
	\label{est:cep'-w}
\end{equation}

Let us consider the Fourier components $u_{k}(t)$ of $u(t)$, and let
us define $E_{k,\ep}(t)$ as in (\ref{defn:ekep}).  Computing the time
derivative as in the strictly hyperbolic case, and using
(\ref{est:cep-c-w}), (\ref{est:cep'-w}), and the fact that
$\cep(t)\geq\omega(\ep)$, we find that
\begin{equation}
	E_{k,\ep}'(t)\leq\gamma_{5}\left(
	\frac{1}{\ep}+\lk\sqrt{\omega(\ep)}\right)E_{k,\ep}(t)
	\quad\quad\forall t\in[0,S].
	\label{est:ekep'}
\end{equation}

Now we choose $\ep$ as a function of $k$. If $\lk<1$ we set
$\ep_{k}=1$, and we have that
$$\gamma_{5}\left(\frac{1}{\ep_{k}}+\lk\sqrt{\omega(\ep_{k})}\right)
\leq\gamma_{5}\left(1+\sqrt{\omega(1)}\right)\leq
\gamma_{6}\varphi(\lk),$$
where we exploited that $\varphi(\sigma)\geq 1$ for all $\sigma\geq 0$.

If $\lk\geq 1$ we consider the function
$h(\sigma)=\sigma\sqrt{\omega(\sigma)}$, which is invertible, and we
set $\ep_{k}:=h^{-1}(1/\lk)$.  By (\ref{defn:Lambda-dg}) we have
that
\begin{equation}
    \lk\sqrt{\omega(\ep_{k})}=\frac{1}{\ep_{k}}\leq\Lambda
    \varphi\left(\frac{1}{h(\ep_{k})}\right)=
    \Lambda\varphi(\lk),
    \label{est:Lambda}
\end{equation}
hence
$$\gamma_{5}\left(\frac{1}{\ep_{k}}+\lk\sqrt{\omega(\ep_{k})}\right)
\leq\gamma_{6}\varphi(\lk).$$

Going back to (\ref{est:ekep'}), in both cases we have that
$E_{k,\ep_{k}}'(t)\leq\gamma_{6}\varphi(\lk)E_{k,\ep_{k}}(t)$ for
every $t\in[0,S]$.  Integrating this differential inequality we obtain
that
\begin{equation}
	E_{k,\ep_{k}}(t)\leq E_{k,\ep_{k}}(0)
	\exp\left(\gamma_{6}\varphi(\lk)t\right)
	\quad\quad\forall t\in[0,S].
	\label{est:ekept}
\end{equation}

In order to estimate $E_{k,\ep_{k}}(0)$ we
need an estimate on $c_{\ep_{k}}(0)$.  To this end we first observe
that 
\begin{equation}
	\omega(\ep_{k})\leq\gamma_{3}
	\quad\quad\forall k\in\n.
	\label{est:omega-ek}
\end{equation}

Indeed this estimate is trivial if $\lk<1$, while for $\lk\geq 1$ it follows
from the fact that $h(\ep_{k})=1/\lk\leq 1$.
Thanks to (\ref{defn:cep-dg}), (\ref{est:c-dg}), and
(\ref{est:omega-ek}) we thus obtain that
$$\cep(0)\leq\gamma_{3}+\m{u_{0}}+1,$$
hence
\begin{equation}
	E_{k,\ep_{k}}(0)\leq \left(\gamma_{3}+\m{u_{0}}+1\right)
	\left(|u_{1k}|^{2}+\lk^{2}|u_{0k}|^{2}\right).
	\label{est:ekepk0}
\end{equation}

Moreover from (\ref{th:omega-3}) and (\ref{est:Lambda}) it follows
that 
\begin{equation}
	\max\left\{1,\frac{1}{\omega(\ep_{k})}\right\}\leq
	1+\frac{1}{\omega(\ep_{k})}\leq
	\gamma_{2}\left(1+\frac{1}{\ep_{k}}\right)\leq
	2\gamma_{2}(1+\Lambda\varphi(\lk)),
	\label{est:oek}
\end{equation}
both for $\lk<1$ and for $\lk\geq 1$.  From (\ref{est:ekept}),
(\ref{est:ekepk0}), and (\ref{est:oek}) it follows that
\begin{eqnarray}
	\lefteqn{\hspace{-5em}
	|u_{k}'(t)|^{2}+\lk^{2}|u_{k}(t)|^{2}\ \leq\ 
	\max\left\{1,\frac{1}{\omega(\ep_{k})}\right\}E_{k,\ep_{k}}(t)}
	\nonumber \\
	& \leq & \gamma_{4}(1+\Lambda\varphi(\lk))
	\left(|u_{1k}|^{2}+\lk^{2}|u_{0k}|^{2}\right)
	\exp\left(\gamma_{6}\varphi(\lk)t\right).
	\label{est:uk-t}
\end{eqnarray}

We are now ready to prove (\ref{th:apriori-est}). We consider the
inequality
\begin{equation}
	1+\Lambda x\leq\left(1+\frac{\Lambda}{M}\right)\exp\left(
	Mx\right)
	\quad\quad\forall\Lambda\geq 0,\ \forall x\geq 0,\ \forall M>0,
	\label{inequality}
\end{equation}
and we apply it with $x=\varphi(\lk)$, $M=r_{0}/2$.  Combined with
(\ref{est:uk-t}) and the fact that $\gamma_{6}T\leq r_{0}/2$ we obtain
that 
$$|u_{k}'(t)|^{2}+\lk^{2}|u_{k}(t)|^{2}\leq
\gamma_{4}\left(1+\frac{2\Lambda}{r_{0}}\right)
\left(|u_{1k}|^{2}+\lk^{2}|u_{0k}|^{2}\right)
\exp\left(r_{0}\varphi(\lk)\right).$$

Multiplying by $\lk$, and summing over all $k$'s, we obtain the same
contradiction as in the strictly hyperbolic case.

The proof of (\ref{th:apriori-reg}) is quite similar. We fix
$\tau\in(0,T]$ and we apply inequality (\ref{inequality}) with
$x=\varphi(\lk)$, $M=\gamma_{6}\tau$. Combined with (\ref{est:uk-t})
we obtain that
\begin{eqnarray*}
	\lefteqn{\hspace{-5em}
	\exp(\strut(r_{0}-R\tau)\varphi(\lk))
	\sup_{t\in[0,\tau]}\left(|u_{k}'(t)|^{2}+\lk^{2}|u_{k}(t)|^{2}\right)
	\ \leq} \\
	 & \leq & \gamma_{4}
	 \left(1+\frac{\Lambda}{\gamma_{6}\tau}\right)
	 \left(|u_{1k}|^{2}+\lk^{2}|u_{0k}|^{2}\right)
	 \exp\left(r_{0}\varphi(\lk)\right)
\end{eqnarray*}
for every $t\in[0,\tau]$.  Multiplying by $\lk$, and summing over all
$k$'s, we conclude as in the strictly hyperbolic case.

%

\label{NumeroPagine}

\end{document}